\documentclass[oneside,leqno,10pt]{article}
\usepackage{amssymb,amsmath,latexsym,amscd}

\usepackage{setspace}
\usepackage[nottoc]{tocbibind}

\newcommand{\diag}{\operatorname{diag}}
 
\def\ga{\mathfrak{a}}

\def\ge{\mathfrak{e}}
\def\gf{\mathfrak{f}}
\def\gg{\mathfrak{g}}
\def\gh{\mathfrak{h}}

\def\gk{\mathfrak{k}}
\def\gl{\mathfrak{l}}
\def\gm{\mathfrak{m}}
\def\gn{\mathfrak{n}}
\def\go{\mathfrak{o}}
\def\gp{\mathfrak{p}}

\def\gs{\mathfrak{s}}
\def\gt{\mathfrak{t}}
\def\gu{\mathfrak{u}}
\def\gv{\mathfrak{v}}

\def\gz{\mathfrak{z}}

\def\C{\mathbb{C}}

\def\R{\mathbb{R}}

\def\Z{\mathbb{Z}}

\def\cC{\mathcal{C}}

\def\cH{\mathcal{H}}

\def\cL{\mathcal{L}}

\def\cO{\mathcal{O}}

\def\Ad{{\rm Ad}}
\def\ad{{\rm ad}\,}
\def\Pf{{\rm Pf}}
\def\rank{{\rm rank}\,}
\def\Ind{{\rm Ind\,}}

\def\tr{{\rm trace\,}}

\def\Det{\rm Det}

\def\diag{{\rm diag}}

\newtheorem{theorem}[equation]{Theorem}

\newtheorem{lemma}[equation]{Lemma}

\newtheorem{corollary}[equation]{Corollary}

\newtheorem{proposition}[equation]{Proposition}

\newtheorem{definition}[equation]{Definition}

\def\sideremark#1{\ifvmode\leavevmode\fi\vadjust{\vbox to0pt{\vss
 \hbox to 0pt{\hskip\hsize\hskip1em
\vbox{\hsize2cm\tiny\raggedright\pretolerance10000 
 \noindent #1\hfill}\hss}\vbox to8pt{\vfil}\vss}}} 

\usepackage{tocloft}

\title{Solvability, Structure and Analysis for Minimal Parabolic Subgroups}

\date{22 January 2017}

\author{Joseph A. Wolf\,\footnote{Research partially supported by a grant 
from the Simons Foundation}}
 
\begin{document}
\maketitle

\begin{abstract}
We examine the structure of the Levi component $MA$ in a minimal parabolic
subgroup $P = MAN$ of a real reductive Lie group $G$ and work out the cases
where $M$ is metabelian, equivalently where $\gp$ is solvable. When 
$G$ is a linear group we verify that $\gp$ is solvable if and only if $M$ 
is commutative.  In the general case $M$ is abelian modulo the center $Z_G$\,,
we indicate the exact structure of $M$ and $P$, and we work out the precise
Plancherel Theorem and Fourier Inversion Formulae.  This lays the groundwork
for comparing tempered representations of $G$ with those induced from
generic representations of $P$.
\end{abstract}

\maketitle

\begin{quote}
{\footnotesize
\begin{spacing}{1.4}
\tableofcontents
\end{spacing}
}
\end{quote}

\section{Introduction}\label{sec1}
\setcounter{equation}{0}

Let $G$ be a real reductive Lie group and $P = MAN$ a minimal parabolic
subgroup.  Later we will be more precise about conditions on the structure
of $G$, but first we recall the {\sl unitary principal series} representations 
of $G$.  They
are the induced representations $\pi_{\chi,\nu,\sigma} =
\text{\rm Ind}_P^G(\eta_{\chi,\nu,\sigma})$ defined as follows.
First, $\nu$ is the highest weight of an irreducible representation 
$\eta_\nu$ of 
the identity component $M^0$ of $M$.  Second, $\chi$ is an irreducible
representation of the $M$--centralizer $Z_M(M^0)$ that agrees with
$\eta_\nu$ on the center $Z_{M^0} = Z_M(M^0) \cap M^0$ of $M^0$.  Third,
$\sigma$ is a real linear functional on the Lie algebra $\ga$ of $A$, in other
words $e^{i\sigma}$ is a unitary character on $A$.  Write
$\eta_{\chi,\nu}$ for the representation 
$\chi \otimes \eta_\nu$ of $M$, and let
$\eta_{\chi,\nu,\sigma}$ denote the representation
$man \mapsto e^{i\sigma}(a)\eta_{\chi,\nu}(m)$ of $P$.  This data
defines the principal series representation $\pi_{\chi,\nu,\sigma} =
\text{\rm Ind}_P^G(\eta_{\chi,\nu,\sigma})$ of $G$.

Now consider a variation in which an irreducible unitary representation of
$N$ is incorporated.  A few years ago, we described Plancherel almost all 
of the unitary dual $\widehat{N}$ in terms of strongly orthogonal roots
(\cite{W2013}, \cite{W2014}).  Using those ``stepwise square integrable''
representations $\pi_\lambda$ of $N$ we arrive at representations 
$\eta_{\chi,\nu,\sigma,\lambda}$ of $P$ and 
$\pi_{\chi,\nu,\sigma,\lambda} = 
\text{\rm Ind}_P^G(\eta_{\chi,\nu,\sigma,\lambda})$ of $G$.
The representations $\pi_{\chi,\nu,\sigma,\lambda}$ have not yet been studied,
at least in terms of their relation to tempered representations and
harmonic analysis on $G$.  In this paper we lay some of the groundwork for
that study. 

Clearly this all this is much simpler when $M$ is commutative modulo the 
center $Z_G$.  Then there is a better chance of finding a clear relation between the
$\pi_{\chi,\nu,\sigma,\lambda}$ and the tempered representation theory of $G$.
In this paper we see just when $M$ is commutative mod $Z_G$.  That
turns out to be equivalent to solvability of $P$, and leads to a straightforward
construction both of the Plancherel Formula and the Fourier Inversion 
Formula for $P$ and of the principal series representations of $G$.
It would also be interesting to see whether solvability of $P$ simplifies the
operator--theoretic formulation \cite{BB2016} of
stepwise square integrability.

It will be obvious to the reader that if any parabolic subgroup of a real
Lie group has commutative Levi component, then that parabolic is a minimal
parabolic.  For this reason we only deal with minimal parabolics.

The concept of ``stepwise square integrable'' representation is basic to
this note and to many of the references after 2012.  It came out of 
conversations with Maria Laura Barberis concerning application of 
square integrability \cite{MW1973} to her work with Isabel Dotti on 
abelian complex structures.  The first developments were \cite{W2013}
and \cite{W2014}, and we follow the notation in those papers.

\section{Lie Algebra Structure}\label{sec2}
\setcounter{equation}{0}

Let $\gg$ be a real reductive Lie algebra.  In other words
$\gg = \gg' \oplus \gz$ where $\gg' = [\gg , \gg]$ is semisimple
and $\gz$ is the center of $\gg$\,.  As usual, $\gg_{_\C}$
denotes the complexification of $\gg$\,,
so $\gg_{_\C} = \gg_{_\C}' \oplus \gz_{_\C}$\,,
direct sum of the respective complexifications of $\gg'$ and
$\gz$\,.  Choose a Cartan involution $\theta$ of $\gg$ and
decompose $\gg = \gk + \gs$ into $(\pm 1)$--eigenspaces of $\theta$.
Fix a maximal abelian subspace $\ga \subset \gs$ and let $\gm$
denote the $\gk$--centralizer of $\ga$.  Let $\gt$ be a Cartan subalgebra
of $\gm$, so $\gh := \gt + \ga$ is a ``maximally split'' Cartan
subalgebra of $\gg$.

We denote root systems by $\Delta(\gg_{_\C},\gh_{_\C})$, 
$\Delta(\gm_{_\C},\gt_{_\C})$ and $\Delta(\gg,\ga)$.  Choose consistent
positive root subsystems $\Delta^+(\gg_{_\C},\gh_{_\C})$,
$\Delta^+(\gm_{_\C},\gt_{_\C})$ and $\Delta^+(\gg,\ga)$.  In other
words, if $\alpha \in \Delta(\gg_{_\C},\gh_{_\C})$ then
$\alpha \in \Delta^+(\gg_{_\C},\gh_{_\C})$ if and only if
either (i) $\alpha|_\ga \ne 0$ and $\alpha|_\ga \in \Delta^+(\gg,\ga)$,
or (ii) $\alpha|_\ga = 0$ and 
$\alpha|_{\gt_{_\C}} \in \Delta^+(\gm_{_\C},\gt_{_\C})$.  We write
$\Psi(\gg_{_\C},\gh_{_\C})$, $\Psi(\gm_{_\C},\gt_{_\C})$ and $\Psi(\gg,\ga)$
for the corresponding simple root systems.  Note that
\begin{equation}\label{simple-simple}
\begin{aligned}
&\Psi(\gg,\ga) = \{\psi|_\ga \mid \psi \in \Psi(\gg_{_\C},\gh_{_\C})
	\text{ and } \psi|_\ga \ne 0 \} \text{ and } \\
&\Psi(\gm_{_\C},\gt_{_\C}) = \{\psi|_{\gt_{_\C}} \mid
	\psi \in \Psi(\gg_{_\C},\gh_{_\C}) \text{ and } 
	\psi|_\ga = 0\}. 
\end{aligned}
\end{equation}
If $\gg$ is the underlying structure of a complex simple Lie algebra $\gl$ 
then $\gg_{_\C} \cong \gl \oplus \overline{\gl}$ and 
$\Psi(\gg_{_\C},\gh_{_\C})$ is the union of the simple root systems of the 
two summands.  In that case $\Psi(\gg,\ga)$ looks like the simple root
system of $\gl$, but with every root of multiplicity $2$.

On the group level, the centralizer $Z_G(\ga) = Z_G(A) = M \times A$ where
$A = \exp(\ga)$ and $M = Z_K(A)$. 

\begin{lemma}\label{gm-abelian-criterion}
The following conditions are equivalent:
$(1)$ if $\alpha \in \Psi(\gg_{_\C},\gh_{_\C})$ then $\alpha|_\ga \ne 0$,
$(2)$ the Lie algebra $\gm$ of $M$ is abelian, and
$(3)$ $\gm$ is solvable.
\end{lemma} 

\noindent {\bf Proof.}  If (1) fails we have
$\psi \in \Psi(\gg_{_\C},\gh_{_\C})$ such that
$\psi|_\ga = 0$, so $0 \ne \psi|_{\gt_{_\C}} \in \Psi(\gm_{_\C},\gt_{_\C})$.
Then $\gm_{_\C}$ contains the simple Lie algebra
with simple root $\psi|_{\gt_{_\C}}$, and (3) fails.
If (3) fails then (2) fails.  
Finally suppose that (2) fails.  Since $\gm_{_\C} = \gt_{_\C} + \sum_{\alpha
\in \Delta(\gm_{_\C},\gt_{_\C})} \gm_\alpha$\,, the root system
$\Delta(\gm_{_\C},\gt_{_\C})$ is not empty.   In particular
$\Psi(\gm_{_\C},\gt_{_\C}) \ne \emptyset$, and (\ref{simple-simple})
provides $\psi \in \Psi(\gg_{_\C},\gh_{_\C})$ with $\psi|_\ga = 0$, so
(1) fails.  \hfill $\square$

Decompose the derived algebra as a direct sum of simple ideals, 
$\gg' = \bigoplus \gg_i$\,.  Then the minimal parabolic subalgebras
$\gp = \gm + \ga + \gn$ 
of $\gg$ decompose as direct sums $\gp = \gz \oplus \bigoplus \gp_i$
where $\gp_i = \gm_i + \ga_i + \gn_i$ is a minimal parabolic 
subalgebra of $\gg_i$\,.  Thus $\gm$ is abelian (resp. solvable) if and
only if each of the $\gm_i$ is abelian (resp. solvable).  The classification
of real reductive Lie algebras $\gg$ with $\gm$ abelian (resp. solvable)
is thus reduced to the case where $\gg$ is simple.  This includes the case
where $\gg$ is the underlying real structure of a complex simple Lie algebra.

The Satake diagram for $\Psi(\gg_{_\C},\gh_{_\C})$ is the
Dynkin diagram, using the arrow convention rather than the black dot
convention, with the following modifications.  If $\psi'$ and $\psi''$
have the same nonzero restriction to $\ga$ then the corresponding
nodes on the diagram are joined by a two--headed arrow.  In the case where
$\gg$ is complex this joins two roots that are complex conjugates of each
other.  If $\psi|_\ga = 0$ 
then the corresponding node on the diagram is changed from a circle to a 
black dot.  Condition (1) of Lemma \ref{gm-abelian-criterion} says that
the Satake diagram for $\Psi(\gg_{_\C},\gh_{_\C})$ has no black dots.
Combining the classification with Lemma \ref{gm-abelian-criterion} we
arrive at

\begin{theorem}\label{s-diagram-list}
Let $\gg$ be a simple real Lie algebra and let $\gp = \gm + \ga + \gn$
be a minimal parabolic subalgebra.  Let $\gt$ be a Cartan subalgebra of
$\gm$, so $\gh := \gt + \ga$ is a maximally split Cartan subalgebra of $\gg$.
Then the following conditions are equivalent$:$ {\rm (i)} $\gm$ is abelian,
{\rm (ii)} $\gm$ is solvable, {\rm (iii)} $\ga$ contains a regular element
of $\gg$, {\rm (iv)} $\gg_{_\C}$ has a Borel subalgebra stable under
complex conjugation of $\gg_{_\C}$ over $\gg$, and {\rm (v)} $\gg$ 
appears on the following list.
\begin{itemize}

\item[$1.$] Cases $\gh = \ga$ and $\gm = 0$ {\rm (called 
split real forms or Cartan normal forms):} 
$\gg$ is one of $\gs\gl(\ell +1;\R)$, $\gs\go(\ell,\ell +1)$,
$\gs\gp(\ell;\R)$, $\gs\go(\ell,\ell)$, $\gg_{2,A_1A_1}$,
$\gf_{4,A_1C_3}$, $\ge_{6,C_4}$, $\ge_{7,A_7}$ or $\ge_{8,D_8}$\,.
In this case, since $\gh = \ga$, the restricted roots all have  
multiplicity $1$.

\item[$2.$] Cases where $\gg$ is the underlying real structure of a 
complex simple Lie algebra$:$ $\gg$ is one of $\gs\gl(\ell +1;\C)$,
$\gs\go(2\ell +1;\C)$, $\gs\gp(\ell;\C)$, $\gs\go(2\ell;\C)$,
$\gg_2$, $\gf_4$, $\ge_6$, $\ge_7$ or $\ge_8$\,.  
Here $\gh = i\ga + \ga$ and $\gm = i\ga$, and
the restricted roots all have multiplicity $2$.

\item[$3.$] Four remaining cases$:$ 
\begin{itemize}
\item[$(3a)$]
$\gg = \gs \gu (\ell, \ell +1)$ with Satake diagram
\setlength{\unitlength}{.5 mm}
\begin{picture}(78,8)
\put(5,-6){\circle{2}}
\put(5,8){\circle{2}}
\put(4,0){$\updownarrow$}
\put(6,-6){\line(1,0){13}}
\put(6,8){\line(1,0){13}}
\put(22,-6){\circle*{1}}
\put(25,-6){\circle*{1}}
\put(28,-6){\circle*{1}}
\put(22,8){\circle*{1}}
\put(25,8){\circle*{1}}
\put(28,8){\circle*{1}}
\put(31,-6){\line(1,0){13}}
\put(31,8){\line(1,0){13}}
\put(45,-6){\circle{2}}
\put(45,8){\circle{2}}
\put(44,0){$\updownarrow$}
\put(46,-6){\line(1,0){6}}
\put(46,8){\line(1,0){6}}
\put(55,-6){\circle*{1}}
\put(58,-6){\circle*{1}}
\put(61,-6){\circle*{1}}
\put(55,6){\circle*{1}}
\put(58,8){\circle*{1}}
\put(61,8){\circle*{1}}
\put(64,-6){\line(1,0){6}}
\put(64,8){\line(1,0){6}}
\put(71,-6){\circle{2}}
\put(71,8){\circle{2}}
\put(70,0){$\updownarrow$}
\put(72,-6){\line(1,0){3}}
\put(72,8){\line(1,0){3}}
\put(75,-6){\line(0,1){14}}
\end{picture}.
In this case $\Delta(\gg, \ga)$ is of type \vskip 7pt $B_\ell$\,,
the long indivisible roots have multiplicity $2$, the short indivisible 
roots have multiplicity $1$, and the divisible roots also have
multiplicity $1$.
\item[$(3b)$]
$\gg = \gs \gu (\ell, \ell)$ with Satake diagram 
\setlength{\unitlength}{.5 mm}
\begin{picture}(90,16)
\put(5,-6){\circle{2}}
\put(5,8){\circle{2}}
\put(4,0){$\updownarrow$}
\put(6,-6){\line(1,0){13}}
\put(6,8){\line(1,0){13}}
\put(22,-6){\circle*{1}}
\put(25,-6){\circle*{1}}
\put(28,-6){\circle*{1}}
\put(22,8){\circle*{1}}
\put(25,8){\circle*{1}}
\put(28,8){\circle*{1}}
\put(31,-6){\line(1,0){13}}
\put(31,8){\line(1,0){13}}
\put(45,-6){\circle{2}}
\put(45,8){\circle{2}}
\put(44,0){$\updownarrow$}
\put(46,-6){\line(1,0){6}}
\put(46,8){\line(1,0){6}}
\put(55,-6){\circle*{1}}
\put(58,-6){\circle*{1}}
\put(61,-6){\circle*{1}}
\put(55,8){\circle*{1}}
\put(58,8){\circle*{1}}
\put(61,8){\circle*{1}}
\put(64,-6){\line(1,0){6}}
\put(64,8){\line(1,0){6}}
\put(71,-6){\circle{2}}
\put(71,8){\circle{2}}
\put(70,0){$\updownarrow$}
\put(72,8){\line(2,-1){13}}
\put(86,0.7){\circle{2}}
\put(72,-6){\line(2,1){13}}
\end{picture}. 
In this case $\Delta(\gg,\ga)$ is of type \vskip 7 pt $C_\ell$\,, 
the long restricted 
roots have multiplicity $1$, and the short ones have multiplicity $2$.
\item[$(3c)$]
$\gg = \gs\go(\ell -1, \ell +1)$ with Satake diagram
\setlength{\unitlength}{.5 mm}
\begin{picture}(80,16)
\put(5,1){\circle{2}}
\put(6,1){\line(1,0){13}}
\put(20,1){\circle{2}}
\put(21,1){\line(1,0){13}}
\put(37,1){\circle*{1}}
\put(40,1){\circle*{1}}
\put(43,1){\circle*{1}}
\put(46,1){\line(1,0){13}}
\put(60,1){\circle{2}}
\put(61,0.7){\line(2,-1){13}}
\put(75,-6){\circle{2}}
\put(61,0.7){\line(2,1){13}}
\put(75,8){\circle{2}}
\put(74,0){$\updownarrow$}
\end{picture}.  
In this case $\Delta(\gg,\ga)$ is of type \vskip 7pt $B_{\ell - 1}$\,, the long
restricted roots have multiplicity $1$, and the short ones have
multiplicity $2$.
\item[$(3d)$]
$\gg = \ge_{6,A_1A_5}$ with Satake diagram
\setlength{\unitlength}{.5 mm}
\begin{picture}(60,16)
\put(5,1){\circle{2}}
\put(6,1){\line(1,0){13}}
\put(20,1){\circle{2}}
\put(21,0.7){\line(2,-1){13}}
\put(35,-6){\circle{2}}
\put(21,0.7){\line(2,1){13}}
\put(35,8){\circle{2}}
\put(34,0){$\updownarrow$}
\put(36,-6){\line(1,0){13}}
\put(36,8){\line(1,0){13}}
\put(50,-6){\circle{2}}
\put(50,8){\circle{2}}
\put(49,0){$\updownarrow$}
\end{picture}.
In this case $\Delta(\gg,\ga)$ is of type $F_4$\,, the long \vskip 7pt 
restricted roots have multiplicity $1$, and the short ones 
have multiplicity $2$.
\end{itemize}
\end{itemize}
$($These real Lie algebras are often called the {\em Steinberg normal forms}
of their complexifications.$)$
\end{theorem}

\begin{corollary}\label{divisibility}
Let $\gg$ be a simple real Lie algebra and let $\gp = \gm + \ga + \gn$
be a minimal parabolic subalgebra with $\gm$ abelian.  Let $\gt$ be a Cartan 
subalgebra of $\gm$, so $\gh := \gt + \ga$ is a maximally split Cartan 
subalgebra of $\gg$.

If $\gg \ne \gs\gu(\ell,\ell +1)$ then
$\Delta(\gg,\ga)$ is nonmultipliable, in other words if 
$\alpha \in \Delta(\gg,\ga)$ then
$2\alpha \notin \Delta(\gg,\ga)$.

If $\gg = \gs\gu(\ell,\ell +1)$ 
\begin{picture}(80,16)
\put(5,-1){\circle{2}}
\put(5,13){\circle{2}}
\put(4,4){$\updownarrow$}
\put(6,-1){\line(1,0){13}}
\put(6,13){\line(1,0){13}}
\put(22,-1){\circle*{1}}
\put(25,-1){\circle*{1}}
\put(28,-1){\circle*{1}}
\put(22,13){\circle*{1}}
\put(25,13){\circle*{1}}
\put(28,13){\circle*{1}}
\put(31,-1){\line(1,0){13}}
\put(31,13){\line(1,0){13}}
\put(45,-1){\circle{2}}
\put(45,13){\circle{2}}
\put(44,4){$\updownarrow$}
\put(46,-1){\line(1,0){6}}
\put(46,13){\line(1,0){6}}
\put(55,-1){\circle*{1}}
\put(58,-1){\circle*{1}}
\put(61,-1){\circle*{1}}
\put(55,13){\circle*{1}}
\put(58,13){\circle*{1}}
\put(61,13){\circle*{1}}
\put(64,-1){\line(1,0){6}}
\put(64,13){\line(1,0){6}}
\put(71,-1){\circle{2}}
\put(71,13){\circle{2}}
\put(70,4){$\updownarrow$}
\put(72,-1){\line(1,0){3}}
\put(72,13){\line(1,0){3}}
\put(75,-1){\line(0,1){14}}
\end{picture},
let $\alpha_1, \dots ,\alpha_{2\ell}$ be the simple roots of
$\Delta^+(\gg_{_\C},\gh_{_\C})$ in the usual order and 
$\psi_i = \alpha_i|_{\ga}$\,.  Then
the multipliable roots in $\Delta^+(\gg,\ga)$ are just the
$\frac{1}{2}\beta_u = (\psi_u + \dots \psi_\ell) = 
(\alpha_u + \dots + \alpha_\ell)|_\ga$ for $1 \leqq u \leqq \ell$; there 
$\beta_u = 2(\psi_u + \dots + \psi_\ell) = 
(\alpha_u + \dots + \alpha_{2\ell -u +1})|_\ga$\,.
\end{corollary}

\section{Structure of the Minimal Parabolic Subgroup} \label{sec3}
\setcounter{equation}{0}
As in Section \ref{sec2}, $\gg$ is a real simple Lie algebra, $\gp =
\gm + \ga + \gn$ is a minimal parabolic subalgebra, and we assume that
$\gm$ is abelian.  $G$ is a connected Lie group with Lie algebra $\gg$
and $P = MAN$ is the minimal parabolic subgroup with Lie algebra $\gp$.
In this section we work out the detailed structure of $M$ and $P$, essentially
by adapting results of K. D. Johnson (\cite{J1987}, \cite{J2004}).
We move the commutativity criteria of 
Theorem \ref{s-diagram-list} from $\gm$ to $M$ when $G$ is linear
and describe the metabelian (in fact abelian mod $Z_G$) 
structure of $M$ for $G$ in general.  

We use the following notation. $\widetilde{G}$ is the connected simply connected
Lie group with Lie algebra $\gg$, $G_{_\C}$ is the connected
simply connected complex Lie group with Lie algebra $\gg_{_\C}$, $G'$ 
is the analytic subgroup of $G_{_\C}$ with Lie algebra $\gg$, and
$\overline{G}$ and $\overline{G_{_\C}}$ are the adjoint groups of
$G$ and $G_{_\C}$.  $G'$ is called the {\em algebraically
simply connected} group for $\gg$.  We write $\widetilde{P} = 
\widetilde{M}\widetilde{A}\widetilde{N}\subset\widetilde{G}$, 
$P' = M'A'N'\subset G'$, and
$\overline{P} = \overline{M}\overline{A}\overline{N}\subset\overline{G}$
for the corresponding minimal parabolic subgroups, aligned so that
$\widetilde{G} \to G'$ maps $\widetilde{M} \to M'$, $\widetilde{A}\cong A'$
and $\widetilde{N} \cong N'$\,, $G' \to \overline{G}$ maps
$M' \to \overline{M}$, $A' \cong \overline{A}$ and $N' \cong \overline{N}$\,,
and $G \to \overline{G}$ maps $M \to \overline{M}$, $A \cong \overline{A}$
and $N \cong \overline{N}$\,.

\subsection{The Linear Case}\label{ssec3a}
In order to discuss $M/M^0$ we need the following concept from \cite{J1987}.
\begin{equation}
\begin{aligned}
r = r(\gg) &\text{ is the
number of white dots in the Satake diagram of $\gg$ } \\
&\text{ not adjacent to a black dot and not attached to another dot 
	by an arrow.}
\end{aligned}
\end{equation}
However we only need it for the case where $\gm$ is abelian, so there are
no black dots.  In fact, running through the cases of Theorem
\ref{s-diagram-list} we have the following.

\begin{proposition}\label{cases-ell}
Let $G$ be a connected Lie group with Lie algebra $\gg$, let $P = MAN$ 
be a minimal parabolic subgroup of $G$, and assume that $\gm$ is abelian.
Then $r(\gg)$ is given by
\begin{itemize}

\item[$1.$] Cases $\gh = \ga$ and $\gm = 0$:  then $r(\gg) = \dim \gh$,
the rank of $\gg$.

\item[$2.$] Cases where $\gg$ is the underlying real structure of a complex
simple Lie algebra:  then $r(\gg) = 0$.

\item[$3.$] Four remaining cases:  
\begin{itemize}
\item[$(3a)$] Case $\gg = \gs \gu (\ell, \ell +1)$: then $r(\gg) = 0$.
\item[$(3b)$] Case $\gg = \gs \gu (\ell, \ell)$: then $r(\gg) = 1$.
\item[$(3c)$] Case $\gg = \gs\go(\ell -1, \ell +1)$: then $r(\gg) = \ell-2$.
\item[$(3d)$] Case $\gg = \ge_{6,A_1A_5}$: then $r(\gg) = 2$.
\end{itemize}
\end{itemize}
\end{proposition}

The first and second assertions in the following Proposition are mathematical 
folklore; we include their proofs for continuity of exposition.  The third 
part is from \cite{J1987}.  Recall that 
$\gg$ satisfies the conditions of {\rm Theorem \ref{s-diagram-list}}.

\begin{proposition}\label{lin-m}
Let $G$ be a connected Lie group with Lie algebra $\gg$, let $P = MAN$
be a minimal parabolic subgroup of $G$, and assume that $\gm$ is abelian.
Let $G'$ be the algebraically simply connected group for $\gg$.

{\rm (1)} $G$ is linear if and only if $G' \to \overline{G}$ factors into
$G' \to G \to \overline{G}$,  

{\rm (2)} if $G$ is linear with minimal
parabolic $P = MAN$, then $M = F\times M^0$ where 
$F \subset (\exp(i\ga)\cap K)$ is an elementary abelian $2$--group,
and 

{\rm (3)} $M' = F' \times M'^0$ where $F' \cong \Z_2^\ell$, and if
$G$ is linear then $F$ is a quotient of $F'$.
\end{proposition}

\noindent {\bf Proof.} If $G$ is linear it is contained in its
complexification, which is covered by $G_{_\C}$\,; (1) follows.

Let $G$ be linear.  The Cartan involution $\theta$ 
of $G$ extends to is complexification $G_c$
and defines the compact real forms $\gg_u := \gk + i\gs$
and $G_u$ of $\gg_{_\C}$ and $G_c$\,.  Here the complexification
$M_cA_c$ of $MA$ is the centralizer of $\ga$ in $G_c$\,.
Its maximal compact subgroup is the centralizer of the torus $\exp(i\ga)$
in the compact connected group $G_u$\,, so it is connected.  Thus
$M^0\exp(i\ga)$ is the centralizer of $\exp(i\ga)$ in $G_u$.  Now the
centralizer $M$ of $\ga$ in $K$ is
$(M^0\exp(i\ga))\cap K = M^0(\exp(i\ga)\cap K)$.  Thus
$M = (\exp(i\ga)\cap K) \cdot M^0$.  

By construction, $\theta$ preserves $\exp(i\ga)\cap K$.  If
$x \in (\exp(i\ga)\cap K)$ then $\theta(x) = x$ because $x \in K$
and $\theta(x) = x^{-1}$ because $x \in \exp(\ga_{_\C})$, so $x = x^{-1}$.
Now $\exp(i\ga)\cap K$ is an elementary abelian $2$--group, so 
$\exp(i\ga)\cap K = F \times (\exp(i\ga)\cap M^0)$ for an elementary
abelian $2$--subgroup $F$, and (2) follows.

Statement (3) is Theorem 3.5 in Johnson's paper \cite{J1987}.
\hfill $\square$

Now we have a chacterization of the linear case.

\begin{theorem} \label{any-group}
Let $G$ be a connected real reductive Lie group and $P = MAN$
a minimal parabolic subgroup.  If $G$ is linear, then 
$M$ is abelian if and only if $\gm$ is abelian, and the 
following conditions are equivalent: {\rm (1)} $M$ is abelian, 
{\rm (2)} $M$ is solvable, and {\rm (3)} $P$ is solvable.
\end{theorem}

\noindent {\bf Proof.}  Since $G$ is linear it has form $G'/Z$ where $G'$
is its algebraically simply connected covering group, the analytic subgroup 
of $G_{_\C}$ for $\gg$.  If $M'$ is abelian then its Lie algebra $\gm$ is
abelian.  Conversely suppose that $\gm$ is abelian.  
Then $M^0$ is abelian and Proposition \ref{lin-m}(2) ensures that $M$
is abelian.  We have proved that if $G$ is linear then $M$ is abelian if and 
only if $\gm$ is abelian.

Statements $(1) \Rightarrow (2) \Rightarrow (3)$
are immediate so we need only prove $(3) \Rightarrow (1)$.  If $P$ is
solvable then $M$ is solvable so $\gm$ is abelian by Theorem 
\ref{s-diagram-list}, and $M$ is abelian as proved just above.
\hfill $\square$

\subsection{The Finite Groups $D_n$ and $D_\gg$}\label{ssec3b}
In order to deal with the nonlinear cases we need certain finite
groups that enter into the description of the component groups of 
minimal parabolics. Those are the $D_n$ for $\gg$ classical or of type
$\gg_2$\,, and $D_\gg$ for the other exceptional cases.

Let $\{e_1, \dots , e_n\}$ be an orthonormal basis of $\R^n$.  Consider
the multiplicative subgroup 
\begin{equation}
D_n = \{\pm e_{i_1} \dots e_{i_{2\ell}} \mid 1 \leqq i_1 < \dots
	< i_{2\ell} \leqq n\}
\end{equation}
of decomposable even invertible elements in the Clifford algebra of $\R^n$.  
It is contained in $Spin(n)$ and has order $2^n$, and we need it for 
$\gg$ classical.  Denote $\overline{D_n} = 
\{\diag(\pm 1, \dots , \pm 1) \mid \det \diag(\pm 1, \dots , \pm 1) = 1\}
\cong \Z_2^{n-1}$.
Then $\overline{D_n}$ is the image of $D_n$ in $SO(n)$ under the usual map
(vector representation) $v: Spin(n) \to SO(n)$.  Here $\{\pm 1\}$ is the 
center and also the derived group of $D_n$\,, and 
$\overline{D_n} \cong D_n/\{\pm 1\}$.  Also, $D_n$ is related to
the $2$--tori of Borel and Serre \cite{BS1953} and to $2$--torsion in integral
cohomology \cite{B1961}.  Note that $D_3$ is isomorphic to the
quaternion group $\{\pm 1, \pm i, \pm j, \pm k\}$.

As usual we write $\widehat{D_n}$ for the unitary dual of $D_n$\,.  
It contains the $2^{n-1}$ characters that factor through $\overline{D_n}$\,.
Those are the $1$--dimensional representations
\begin{equation}
\varepsilon_{i_1,\dots,i_{2\ell}}: \diag(a_1,\dots,a_n) \mapsto
	a_{i_1}a_{i_2}\dots a_{i_{2\ell}}.
\end{equation}
There are also representations of degree $> 1$:  
\begin{equation}
\begin{aligned}
&\text{If $n = 2k$ even
let $\sigma_\pm$ denote the restriction of the half--spin representations
from $Spin(n)$ to $D_n$\,.}\\
&\text{If $n = 2k+1$ odd let $\sigma$ denote the restriction of the spin
representation from $Spin(n)$ to $D_n$\,.}
\end{aligned}
\end{equation}
Note that $\deg \sigma_\pm = 2^{k-1}$ for $n=2k$ and 
$\deg \sigma = 2^k$ for $n = 2k+1$.
These representations enumerate $\widehat{D_n}$, as follows.  We will need this
for the Plancherel formula for $P$.
\begin{proposition}\label{dual-dn}
{\rm (\cite[Section 3]{J2004})}
The representations $\sigma_\pm$ of $D_{2k}$ and $\sigma$ of $D_{2k+1}$
are irreducible, and
\begin{equation}
\begin{aligned}
&\text{If $n=2k$ even then } \widehat{D_n} = 
        \{\sigma_+\,, \sigma_-\,, \varepsilon_{i_1,\dots,i_{2\ell}} \mid
        1 \leqq i_1 < \dots < i_{2\ell} \leqq n \text{ and } 
        0 \leqq \ell < k\}. \\
&\text{If $n = 2k+1$ odd then } \widehat{D_n} =
        \{\sigma\,, \varepsilon_{i_1,\dots,i_{2\ell}} \mid
        1 \leqq i_1 < \dots < i_{2\ell} \leqq n \text{ and } 
        0 \leqq \ell \leqq k\}.
\end{aligned}
\end{equation}
\end{proposition}

Now we describe the analogs of the $D_n$ for the split exceptional
Lie algebras of types $\gf_4$\,, $\ge_6$\,, $\ge_7$\,, or $\ge_8$\,.  
Following \cite[Section 8]{J2004} the natural inclusions
$\gf_4 \subset \ge_6 \subset \ge_7 \subset \ge_8$ exponentiate to
inclusions 
$$
\widetilde{Z}_{E_{8,D_8}} \subset \widetilde{F}_{4,C_1C_3} \subset 
	\widetilde{E}_{6,C_4} \subset \widetilde{K}_{E_{7,A_7}} \subset 
	\widetilde{E}_{8,D_8}
$$
where $\widetilde{Z}_{E_{8,D_8}}$ is the center of $\widetilde{E}_{8,D_8}$
and the others are split real simply connected exceptional Lie groups.
For the first inclusion, note that the maximal compact subgroups satisfy
$$
\widetilde{K}_{F_{4,C_1C_3}} \subset \widetilde{K}_{E_{6,C_4}} \subset 
\widetilde{K}_{E_{7,A_7}} \subset \widetilde{K}_{E_{8,D_8}} \text{ given by }
(Sp(1)\times Sp(3)) \subset Sp(4) \subset SU(8) \subset Spin(16)
$$
and $\widetilde{Z}_{E_{8,D_8}} = \{1,e_1\cdot e_2 \cdot \dots \cdot e_{16}\}
\cong \Z_2$ in the spin group using Clifford multiplication.

Let $U_7$ denote the group of permutations of $\{1,2,\dots ,8\}$ generated by
products of $4$ commuting transpositions, e.g. by $\tau := (12)(34)(56)(78)$ 
and its conjugates in the permutation group, viewed as a subgroup 
$\cong \Z_2^3$ of $SU(8)$.  Let $V_7$ be the group 
generated by $\omega_1 := iI$, $\omega_2 := \diag(-1,-1,1,1,1,1,-1,-1)$, 
$\omega_3 := \diag(-1,-1,1,1,-1,-1,1,1)$ and
$\omega_4 := \diag(-1,1,-1,1,-1,1,-1,1)$, viewed as a subgroup
$\cong \Z_4 \times \Z_2^3$ of $SU(8)$.  
Now let $U_6$ denote the subgroup $\{1, (13)(24)(57)(68),\, (15)(26)(37)(48),\,
(17)(28)(37)(48)\} \subset U_7$, so $U_6\cong Z_2^2$\, and let $V_6$ denote
the subgroup of $V_7$ generated by $\omega_1\omega_4$\,, $\omega_2$\, 
and $\omega_3$\,.  Finally let $W_4$ denote the group generated by
$\{\tau\omega_1\,, \omega_2\,, \omega_3\,, \tau\omega_4\}$.

\begin{proposition}\label{e4678} {\rm (\cite[Section 9]{J2004}}
Define $W_6 = U_6V_6 \cup \tau\omega_1 U_6V_6$\,, 
$W_7 = U_7V_7$ and $W_8 = W_7 \cup \tau W_7$\,.  Then

{\rm 1.} $W_4$ is a group of order $2^5$ with $[W_4,W_4] = \{\pm 1\}
\cong \Z_2$\,, $W_4/[W_4,W_4] \cong \Z_2^4$\,, and $W_4$ has center
$Z_{W_4} = \{\pm 1, \pm \omega_2\,, \pm\omega_3\,, \pm\omega_2\omega_3\}$.  
The action of $W_4$ on $\C^8$\,, as a subgroup of $Sp(4)$, breaks into 
the sum of $4$ irreducible inequivalent $2$--dimensional subspaces
$\C^2_j$\,, distinguished by the action of $Z_{W_4}$.  Thus the unitary dual 
$\widehat{W_4} = \{w_{4,1}\,, w_{4,2}\,, w_{4,3}\,, w_{4,4}\,,
\epsilon_1, \dots , \epsilon_{16}\}$ where the $w_{4,j}$ are the
representations on the $\C^2_j$ and the $\epsilon_j$ are the
$(1$--dimensional$)$ representations that
annihilate $[W_4,W_4]$.

{\rm 2.} $W_6$ is a group of order $2^7$ with $[W_6,W_6] = 
\{\pm 1\} \cong \Z_2$ and $W_6/[W_6,W_6] \cong \Z_2^6$\,.
Let $w_6 = w|_{W_6}$ where $w$ is the $($vector$)$ representation of $Sp(4)$
on $\C^8$.  Then $w_6$ is irreducible and $\widehat{W_6} =
\{w_6\,, \varepsilon_1\,, \dots , \varepsilon_{64}\}$ where the      
$\varepsilon_j$  are the $(1$--dimensional$)$ representations that
annihilate $[W_6,W_6]$.

{\rm 3.} $W_7$ is a group of order $2^8$ with derived group
$[W_7,W_7] = \widetilde{Z}_{E_{8,D_8}}$\,, and $W_7/[W_7,W_7] \cong \Z_2^7$\,.
Let $w_7 = w|_{W_7}$ where $w$ is the $($vector$)$ representation of $SU(8)$
on $\C^8$.  Then $w_7$ and $w_7^*$ are inequivalent irreducible
representations of $W_7$\,, and the unitary dual
$\widehat{W_7} = \{w_7\,, w_7^*\,, \varepsilon'_1\,, \dots , \varepsilon'_{128}\}$
where the $\varepsilon'_j$ are the $(1$--dimensional$)$ representations that
annihilate $[W_7,W_7]$.

{\rm 4.} $W_8$ is a group of order $2^9$ with $[W_8,W_8] = 
\widetilde{Z}_{E_{8,D_8}}$\, and $W_8/[W_8,W_8] \cong \Z_2^8$\,.
Let $w_8 = w|_{W_8}$ where $w$ is the $($vector$)$ representation of $Spin(16)$
on $\C^{16}$.  Then $w_8$ is irreducible and $\widehat{W_8} =
\{w_8\,, \varepsilon''_1\,, \dots , \varepsilon''_{256}\}$ where the 
$\varepsilon''_j$  are the $(1$--dimensional$)$ representations that
annihilate $[W_8,W_8]$.
\end{proposition}

\subsection{The General Case}\label{ssec3c}
As before, $G$ is a connected real simple Lie group with minimal
parabolic $P = MAN$ such that $\gm$ is abelian.
Recall that $\widetilde{G}$ is the universal covering group of $G$,
$G'$ is the algebraically simply connected Lie group with Lie
algebra $\gg$,
and $\overline{G}$ is the adjoint group.  Also, $\widetilde{P}$, $P'$,
$P$ and $\overline{P}$ are the respective minimal parabolics.
We specialize the summary section of \cite{J2004} to our setting.
\begin{itemize}
\item[$(a)$] $\widetilde{M} = \widetilde{F} \times \widetilde{M}^0$ where 
$\widetilde{F}$ is discrete, and if $\gg$ is the split real form of 
$\gg_{_\C}$ then $\widetilde{M} = \widetilde{F}$ discrete,
\item[$(b)$] $\widetilde{F}$ is infinite if and only if $G/K$ is a 
tube domain (hermitian symmetric space of tube type),
\item[$(c)$] $r(\gg) = 0 \Leftrightarrow \widetilde{M} \text{ is connected }
\Leftrightarrow M' \text{ is connected},$
\item[$(d)$] if $r(\gg) = 1$ then $G/K$ is a tube domain and 
$\widetilde{F} \cong \Z$,
\item[$(e)$] if $r(\gg) > 1$ and $G/K$ is a tube domain then $G' = Sp(n;\R)$
and $\widetilde{M} = \widetilde{F} \cong \Z_2^{r(\gg) - 1} \times \Z$, and
\item[$(f)$] if $r(\gg) > 1$ and $G/K$ is not a tube domain then $\widetilde{F}$
is a nonabelian group of order $2^{r(\gg) + 1}$.
\end{itemize} 
Now we combine this information with Theorems \ref{s-diagram-list}
and \ref{cases-ell}, as follows.  We use $\gs\gu(1,1) = \gs\gp(1;\R) 
= \gs\gl(2;\R)$, and $\gs\go(2,3) = \gs\gp(2;\R)$.

\begin{proposition}\label{mtilde}
Let $G$ be a connected simple Lie group with Lie algebra $\gg$, let $P = MAN$
be a minimal parabolic subgroup of $G$, and assume that $\gm$ is abelian.
Retain the notation $\widetilde{G}$, $\widetilde{P}$ and
$\widetilde{M} = \widetilde{F} \times \widetilde{M}^0$ as above.
\begin{itemize}

\item[$1.$] Cases $\gh = \ga$, where $G$ is a split real Lie group and $r(\gg)
= \rank \gg$.  Then $\widetilde{M} = \widetilde{F}$. 
\begin{itemize}
\item[$(1a)$] If $\widetilde{M}$ is infinite then
$\gg = \gs\gp(n;\R)$, where $n \geqq 1$ and $r(\gg) = n$. 
In that case $\widetilde{M} \cong \Z_2^{n-1} \times Z$.
\item[$(1b)$]If $\widetilde{M}$ is finite then $\widetilde{F}$
is a nonabelian {\rm (but metabelian)} group of order $2^{r(\gg) + 1}$.
In those cases
\begin{itemize}
\item $\gs\gl(n;\R),\, n = 3 \text{ or } n > 4$: $\widetilde{M} \cong D_n$ from 
	{\rm \cite[Proposition 17.1]{J2004}.}
\item $\gs\go(\ell,\ell+1) \text{ and } \gs\go(\ell,\ell),\, \ell \geqq 3$: 
	$\widetilde{M} \cong D_\ell$ from {\rm \cite[Proposition 17.5]{J2004}.}
\item $\gg_{2,A_1A_1}$: $\widetilde{K} = Sp(1)\times Sp(1)$ and
	$\widetilde{M} \cong D_3$ from {\rm \cite[Proposition 10.4]{J2004}.}
\item $\gf_{4,A_1C_3}$: $\widetilde{K} = Sp(1)\times Sp(3)$ and
	$\widetilde{M} \cong W_4$ {\rm (Proposition \ref{e4678} above)}
	from {\rm \cite[Proposition 9.6]{J2004}.}
\item $\ge_{6,C_4}$: $\widetilde{K} = Sp(4)$ and
	$\widetilde{M} \cong W_6$ {\rm (Proposition \ref{e4678} above)}
        from {\rm \cite[Proposition 9.5]{J2004}.}
\item $\ge_{7,A_7}$: $\widetilde{K} = SU(8)$ and 
	$\widetilde{M} \cong W_7$ {\rm (Proposition \ref{e4678} above)}
        from {\rm \cite[Proposition 9.3]{J2004}.}
\item $\ge_{8,D_8}$: $\widetilde{K} = Spin(16)$ and
	$\widetilde{M} \cong W_8$ {\rm (Proposition \ref{e4678} above)}
        from {\rm \cite[Proposition 9.4]{J2004}.}
\end{itemize}
\end{itemize}

\item[$2.$]  Cases where $\gg$ is the underlying real structure of a 
complex simple Lie algebra.  Then $r(\gg) = 0$, $\widetilde{F} = \{1\}$, 
$\widetilde{G} = G'$ {\rm (so $G$ is linear)}, and $\widetilde{M} = M'$
is a torus group.

\item[$3.$] The four remaining cases$:$ 
\begin{itemize}
\item[$(3a)$] $\gg = \gs \gu (\ell, \ell +1)$: $r(\gg) = 0$, 
$\widetilde{F} = \{1\}$ and $\widetilde{M} \cong U(1)^{\ell -1}\times \R$
	from {\rm \cite[Proposition 17.3]{J2004}.}
\item[$(3b)$] $\gg = \gs \gu (\ell, \ell)$: $r(\gg) = 1$, $\widetilde{F}
	\cong \Z$ and $\widetilde{M} \cong U(1)^{\ell -1}\times \Z$
	from {\rm \cite[Proposition 17.4]{J2004}.}
\item[$(3c)$] $\gg = \gs\go(\ell -1, \ell +1)$, $\ell \ne 1$:  
	$r(\gg) = \ell -2$, and
   \begin{itemize} 
	\item if $\ell \ne 3$ then $\widetilde{F} \cong D_{\ell - 1}$ and
        	$\widetilde{M} \cong D_{\ell - 1} \times \R$
		from {\rm \cite[Proposition 17.5]{J2004},}
	\item if $\ell = 3$ and then $\widetilde{F} \cong \Z$ and
		$\widetilde{M} \cong \Z \times \R$ 
		from {\rm \cite[Proposition 17.6]{J2004}.}
   \end{itemize}
\item[$(3d)$] $\gg = \ge_{6,A_1A_5}$: $r(\gg) = 2$, $\widetilde{F} \cong
	D_3$ and $\widetilde{M} \cong D_3 \times Spin(8)$ from
	{\rm \cite[\S 16]{J2004}} and {\rm \cite[Theorem 3.5]{J1987}}
\end{itemize}
\end{itemize}
\end{proposition}

\section{Stepwise Square Integrable Representations of the Nilradical}
\label{sec4}
\setcounter{equation}{0}

In this section we recall the background for Fourier Inversion on 
connected simply connected nilpotent Lie groups,
and its application to nilradicals of minimal parabolics.
In Section \ref{sec5}  we use it to describe ``generic'' representations 
and the Plancherel Formula on the minimal parabolic.  Then in Section \ref{sec6}
we come to the Fourier Inversion Formula on the parabolic.
In the last section (the Appendix) we will 
run through the cases of Theorem \ref{s-diagram-list}
making explicit the Plancherel Formula and Fourier Inversion Formula
for those nilradicals.

The basic decomposition is
\begin{equation}\label{setup}
\begin{aligned}
N = &L_1L_2\dots L_m \text{ where }\\
 &\text{(a) each factor $L_r$ has unitary representations that are square
integrable modulo its center } Z_r\,, \\
 &\text{(b) each } L_r = L_1L_2\dots L_r \text{ is a normal subgroup of } N
	\text{ with } N_r = N_{r-1}\rtimes L_r \text{ semidirect,} \\
&\text{(c) decompose }\gl_r = \gz_r + \gv_r \text{ and } \gn = \gs + \gv
        \text{ as vector direct sums where } \\
 &\phantom{XXXX}\gs = \oplus\, \gz_r \text{ and } \gv = \oplus\, \gv_r;
    \text{ then } [\gl_r,\gz_s] = 0 \text{ and } [\gl_r,\gl_s] \subset \gv
        \text{ for } r > s\,.
\end{aligned}
\end{equation}

We will need the notation
\begin{equation}\label{c-d}
\begin{aligned}
&\text{(a) }d_r = \tfrac{1}{2}\dim(\gl_r/\gz_r) \text{ so }
        \tfrac{1}{2} \dim(\gn/\gs) = d_1 + \dots + d_m\,,
        \text{ and } c = 2^{d_1 + \dots + d_m} d_1! d_2! \dots d_m!\\
&\text{(b) }b_{\lambda_r}: (x,y) \mapsto \lambda([x,y])
        \text{ viewed as a bilinear form on } \gl_r/\gz_r \\
&\text{(c) }S = Z_1Z_2\dots Z_m = Z_1 \times \dots \times Z_m \text{ where } Z_r
        \text{ is the center of } L_r \\
&\text{(d) }P: \text{ polynomial } P(\lambda) = \Pf(b_{\lambda_1})
        \Pf(b_{\lambda_2})\dots \Pf(b_{\lambda_m}) \text{ on } \gs^* \\
&\text{(e) }\gt^* = \{\lambda \in \gs^* \mid P(\lambda) \ne 0\} \\
&\text{(f) } \pi_\lambda \in \widehat{N} \text{ where } \lambda \in \gt^*:
    \text{ irreducible unitary rep. of } N = L_1L_2\dots L_m\,.
\end{aligned}
\end{equation}

As $\exp : \gn \to N$ is a polynomial diffeomorphism, the Schwartz space
$\cC(N)$ consists of all $C^\infty$ functions $f$ on $N$ such that
$f\cdot\exp \in \cC(\gn)$, the classical Schwartz space of all rapidly 
decreasing $C^\infty$ functions on the real vector space $\gn$.
The general result, which we will specialize, is \cite{W2013}
\begin{theorem}\label{plancherel-general}
Let $N$ be a connected simply connected nilpotent Lie group that
satisfies {\rm (\ref{setup})}.  Then Plancherel measure for $N$ is
concentrated on $\{[\pi_\lambda] \mid \lambda \in \gt^*\}$ where
$[\pi_\lambda]$ denotes the unitary equivalence class of $\pi_\lambda$\,.
If $\lambda \in \gt^*$, and if $u$ and $v$ belong to the
representation space $\cH_{\pi_\lambda}$ of $\pi_\lambda$,  then
the coefficient $f_{u,v}(x) = \langle u, \pi_\nu(x)v\rangle$
satisfies
\begin{equation}
||f_{u,v}||^2_{\cL^2(N/S)} = \frac{||u||^2||v||^2}{|P(\lambda)|}\,.
\end{equation}
The distribution character $\Theta_{\pi_\lambda}$ of $\pi_{\lambda}$ satisfies
\begin{equation}
\Theta_{\pi_\lambda}(f) = c^{-1}|P(\lambda)|^{-1}\int_{\cO(\lambda)}
        \widehat{f_1}(\xi)d\nu_\lambda(\xi) \text{ for } f \in \cC(N)
\end{equation}
where $c$ is given by {\rm (\ref{c-d})(a)},
$\cC(N)$ is the Schwartz space, $f_1$ is the lift
$f_1(\xi) = f(\exp(\xi))$, $\widehat{f_1}$ is its classical Fourier transform,
$\cO(\lambda)$ is the coadjoint orbit $\Ad^*(N)\lambda = \gv^* + \lambda$,
and $d\nu_\lambda$ is the translate of normalized Lebesgue measure from
$\gv^*$ to $\Ad^*(N)\lambda$.  The Plancherel Formula on $N$ is
\begin{equation}
\cL^2(N) = \int_{\gt^*} \cH_{\pi_\lambda} \widehat{\otimes} \cH^*_{\pi_\lambda}
	|P(\lambda)|d\lambda
	\text{ where $\cH_{\pi_\lambda}$ is the representation space of }
	\pi_\lambda
\end{equation}
and the Fourier Inversion Formula is
\begin{equation}
f(x) = c\int_{\gt^*} \Theta_{\pi_\lambda}(r_xf) |P(\lambda)|d\lambda
        \text{ for } f \in \cC(N) \text{ with } 
	c \text{ as in  {\rm (\ref{c-d})(a)}}.
\end{equation}
\end{theorem}
\begin{definition}\label{stepwise2}
{\rm The representations $\pi_\lambda$ of (\ref{c-d}(f)) are the
{\it stepwise square integrable} representations of $N$ relative to
(\ref{setup}).}\hfill $\diamondsuit$
\end{definition}

Nilradicals of minimal parabolics fit this pattern as follows.
Start with the Iwasawa decomposition $G = KAN$.  
Here we use the root order of (\ref{simple-simple}) so $\gn$ is the sum 
of the positive $\ga$--root spaces in $\gg$.  Since $\Delta(\gg,\ga)$ is 
a root system, if $\gamma \in \Delta(\gg_{_\C},\gh_{_\C})$ and 
$\gamma|_\ga \in \Delta^+(\gg,\ga)$ then 
$\gamma \in \Delta^+(\gg_{_\C},\gh_{_\C})$.
Define $\beta_1$ to be the maximal root,
$\beta_{r+1}$ a maximum among the positive roots orthogonal to 
$\{\beta_1, \dots , \beta_r\}$, etc.  This constructs a maximal set
$\{\beta_1, \dots , \beta_m\}$ of strongly orthogonal 
positive restricted roots.  For $1\leqq r \leqq m$ define
\begin{equation}\label{layers}
\begin{aligned}
&\Delta^+_1 = \{\alpha \in \Delta^+(\gg,\ga) \mid \beta_1 - \alpha \in \Delta^+(\gg,\ga)\}
\text{ and }\\
&\Delta^+_{r+1} = \{\alpha \in \Delta^+(\gg,\ga) \setminus (\Delta^+_1 \cup \dots \cup \Delta^+_r)
        \mid \beta_{r+1} - \alpha \in \Delta^+(\gg,\ga)\}.
\end{aligned}
\end{equation}
Then 
\begin{equation}\label{redef-del}
\Delta^+_r\cup \{\beta_r\}
= \{\alpha \in \Delta^+ \mid \alpha \perp \beta_i \text{ for } i < r
\text{ and } \langle \alpha, \beta_r\rangle > 0\}.
\end{equation}  
{\em Note: if $\beta_r$ is divisible then $\frac{1}{2}\beta_r\in\Delta^+_r$\,;
See Corollary \ref{divisibility}.}
Now define
\begin{equation}\label{def-l}
\gl_r = \gg_{\beta_r} + {\sum}_{\Delta^+_r}\, \gg_\alpha
\text{ for } 1\leqq r\leqq m.
\end{equation}
Thus $\gn$ has an increasing foliation by ideals
\begin{equation}\label{def-filtration}
\gn_r = \gl_1 + \gl_2 + \dots + \gl_r \text{ for } 1 \leqq r \leqq m.
\end{equation}
The corresponding group level decomposition
$N = L_1L_2\dots L_m$ and the semidirect product decompositions
$N_r = N_{r-1}\rtimes L_r$ satisfy all the requirements of (\ref{setup}).

\section{Generic Representations of the Parabolic}\label{sec5}
\setcounter{equation}{0}

In this section $G$ is a connected real reductive Lie group, not necessarily
linear, and the minimal parabolic subgroup $P = MAN$ is solvable. In other 
words we are in the setting of
Theorems \ref{any-group} and \ref{plancherel-general}.
Recall $\gt^* = \{\lambda = (\lambda_1 + \dots + \lambda_m) \in \gs^* \mid
\text{ each } \lambda_r \in \gg_{\beta_r} \text{ with } 
\Pf_{\gl_r}(\lambda_r) \ne 0\}$.  For each $\lambda \in \gt^*$ we have the
stepwise square integrable representation $\pi_\lambda$ of $N$.  Now we look
at the corresponding representations of $MAN$.  As before, the 
superscript ${}^0$ denotes identity component.

Theorem \ref{any-group} shows that $\Ad(MA)$ is commutative.  Also, $\Ad(M^0A)$
acts $\C$--irreducibly on each complexified restricted root space
$(\gg_\alpha)_{_\C}$ by \cite[Theorem 8.13.3]{W1966}\,.  But $\Ad(M^0A)$
preserves each $\gg_\alpha$, so it is irreducible there.  Thus, for each
$\alpha \in \Delta(\gg,\ga)$, either $\Ad(M^0)$ is trivial on $\gg_\alpha$
and $\dim_{_\R}\gg_\alpha = 1$, or $\Ad(M^0)$ is nontrivial on 
$\gg_\alpha$ and
$\dim_{_\R}\gg_\alpha = 2$.  In the latter case $\Ad(M^0)|_{\gg_\alpha}$ 
must be
the circle group of all proper rotations of $\gg_\alpha$\,.  A glance at
Theorem \ref{s-diagram-list} makes this more explicit on the $\gg_{\beta_r}$\,.

\begin{lemma}\label{m0-action}
In cases {\rm (1)} and {\rm (3)} of {\rm Theorem \ref{s-diagram-list}},
each $\dim_{_\R} \gg_{\beta_r} = 1$, so $\Ad^*(M^0)$ is trivial
on each $\gg_{\beta_r}$\,.  In case {\rm (2)} each 
$\dim_{_\R} \gg_{\beta_r} = 2$, so $\Ad^*(M^0)$ acts on $\gg_{\beta_r}$
as a circle group $SO(2)$.
\end{lemma}
Recall the notation $\widetilde{M} = \widetilde{F}\times \widetilde{M}^0$ 
from Section \ref{ssec3a} where $p: \widetilde{G} \to G$ is the universal
cover and $\widetilde{P} = \widetilde{M} \widetilde{A} \widetilde{N}$
is $p^{-1}(P)$.  Note $M = F\cdot M^0$ where $F = p(\widetilde{F})$.
Combining \cite[Lemma 3.4]{W2014} with \cite[Proposition 3.6]{W2014}
and specializing to the case of $\gm$ abelian, we have
\begin{lemma}\label{f-action}
When $\gm$ is abelian, $\Ad^*(\widetilde{F})$ acts trivially on $\gs^*$. 
Then in particular, $\Ad^*(F)|_{\gs^*}$ is trivial and 
$\Ad^*(M)|_{\gs^*} = \Ad^*(M^0)|_{\gs^*}$\,.
\end{lemma}
Now combine Lemmas \ref{m0-action} and \ref{f-action}:
\begin{proposition}\label{m-action}
In cases {\rm (1)} and {\rm (3)} of {\rm Theorem \ref{s-diagram-list}},
$\Ad^*(M)$ is trivial on $\gs^*$.  In case {\rm (2)} of 
{\rm Theorem \ref{s-diagram-list}}, $\Ad^*(M)$ acts nontrivially as a circle 
group on each $\gg^*_{\beta_r}$\,, thus acts almost--effectively as a torus
group on $\gs^*$.
\end{proposition}

Fix $\lambda \in \gt^*$.  By Proposition \ref{m-action} its
$\Ad^*(M)$--stabilizer is all of $M$ in cases {\rm (1)} and {\rm (3)} of
Theorem \ref{s-diagram-list}, and in case {\rm (2)} it has form
\begin{equation}\label{m-stab}
M_\diamond = FM^0_\diamond \text{ where } 
M^0_\diamond = \{x \in M^0 \mid \Ad(x)|_{\gs} = 1\}.
\end{equation}
This is independent of the choice of $\lambda \in \gt^*$.  
Thus the kernel of the action of $\Ad(M^0)$ on $\gs^*$ is the codimension $m$
subtorus of $M^0$ with Lie algebra $\gm_\diamond = \sqrt{-1}
\{\xi \in \ga \mid \text{ every } \beta_r(\xi) = 0\}$.

Since $\Ad^*(A)$ acts on $\gg_{\beta_r}$ by positive real scalars, 
given by the real character $e^{\beta_r}$, 
we have a similar result for $A$: the $\Ad^*(A)$--stabilizer
of any $\lambda \in \gt^*$ is
\begin{equation}\label{a-stab}
A_\diamond = \{\exp(\xi) \mid \xi \in \ga \text{ and every }
	\beta_r(\xi) = 0\}.
\end{equation}
Its Lie algebra is $\ga_\diamond = 
\{\xi \in \ga \mid \text{ every }\beta_r(\xi) = 0\}$.
Combining (\ref{m-stab}) and (\ref{a-stab}) we arrive at
\begin{lemma}\label{ma-stabilizer}
The stepwise square integrable representations $\pi_\lambda$ of
$N$ all have the same $MA$--stabilizer $M_\diamond A_\diamond$
on the unitary dual $\widehat{N}$.
\end{lemma}

Specialize \cite[Lemma 3.8]{W2014} and \cite[Lemma 5.4]{W2014} to
$\pi_\lambda$ and $M_\diamond A_\diamond N$.  The Mackey
obstruction (\cite{M1952}, \cite{M1958}, or see \cite{M1955}) vanishes
as in \cite{S1971} and \cite{W1975}.  Now $\pi_\lambda$
extends to an irreducible unitary representation $\widetilde{\pi_\lambda}$
of $M_\diamond A_\diamond N$ on the representation space $\cH_{\pi_\lambda}$
of $\pi_\lambda$\,.  Compare \cite{D1972}.  
Consider the unitarily induced representations
\begin{equation}\label{def-induced}
\pi_{\chi,\alpha,\lambda} := \Ind_{M_\diamond A_\diamond N}^{MAN}\,
	(\chi \otimes e^{i\alpha}\otimes\widetilde{\pi_\lambda}) \text{ for }
	\chi \in \widehat{M_\diamond}\,,\, \alpha \in \ga_\diamond^* 
	\text{ and }\lambda \in \gt^*.
\end{equation}
Note that $\chi$ is a (finite--dimensional) unitary representation of the
metabelian group $M_\diamond$ and that the representation
space of $\chi \otimes e^{i\alpha}\otimes \widetilde{\pi_\lambda}$ is 
$\cH_\chi \widehat{\otimes}\, \C \widehat{\otimes} \cH_{\pi_\lambda}$.
So the representation space of $\pi_{\chi,\alpha,\lambda}$ is
\begin{equation}\label{repspace}
\begin{aligned}
\cH_{\pi_{\chi,\alpha,\lambda}} = 
&\{\text{$L^2$ functions }
	f:MAN \to \cH_\chi \widehat{\otimes} \cH_{\pi_\lambda} \mid \\
	&f(xman) = \delta(a)^{-1/2}e^{-i\alpha(\log a)}(\chi(m)^{-1}\otimes
	\widetilde{\pi_\lambda}(man)^{-1}) (f(x))\,, 
	x \in MAN, man \in M_\diamond A_\diamond N\}.
\end{aligned}
\end{equation}
Unitarity of $\pi_{\chi,\alpha,\lambda}$ requires the $\delta(a)^{-1/2}$ term, 
as we'll see when we discuss Dixmier--Puk\' anszky operators.

\begin{definition}\label{def-generic}
{\rm The $\pi_{\chi,\alpha,\lambda}$ of
(\ref{def-induced}) are the {\em generic} irreducible unitary representations
of $MAN$.}
\end{definition}
Now, specializing \cite[Theorem 5.12]{W2014},
\begin{theorem}\label{sup-planch}
Let $G$ be a connected real reductive Lie group and $P=MAN$ a minimal
parabolic subgroup.  Suppose that $P$ is solvable.  Then the
Plancherel measure for $MAN$ is concentrated on the set of all generic
unitary representation classes 
$[\pi_{\chi,\alpha,\lambda}] \in \widehat{MAN}$, and
$$
\cL^2(MAN) = \int_{\chi \in \widehat{M_\diamond}}
	\left ( \int_{\alpha \in \ga_\diamond^*} \left (
	\int_{\lambda \in \gt^*} \cH_{\pi_{\chi,\alpha,\lambda}} 
	\widehat{\otimes} \cH_{\pi_{\chi,\alpha,\lambda}}^* |P(\lambda)|
	d\lambda\right ) d\alpha \right )\deg(\chi)d\chi.
$$
\end{theorem}

\section{Fourier Inversion on the Parabolic}\label{sec6}
\setcounter{equation}{0}
In this section, as before, $G$ is a connected\footnote{These results extend
{\it mutatis mutandis} to all real reductive Lie groups $G$ such
that (a) the minimal parabolic subgroup of $G^0$ is solvable, 
(b) if $g \in G$ then $\Ad(g)$ is an inner automorphisms of $\gg_{_\C}$
and (c) $G$ has a closed normal abelian subgroup $U$ such that (c1) $U$ 
centralizes the identity component $G^0$, (c2) $UG^0$ has finite index in $G$
and (c3) $U\cap G^0$ is co--compact in the center of $G^0$.  The extension
is relatively straightforward using the methods of \cite{W1974} as described 
in \cite[Introduction]{W1974} and \cite[Section 1]{W1974}.  For continuity of 
exposition we the leave details on that to the interested reader.}
real reductive Lie group whose
minimal parabolic subgroup $P = MAN$ is solvable.  We work out an explicit 
Fourier Inversion Formula for $MAN$.  It uses the
generic representations of Definition \ref{def-generic} and an operator to
compensate non--unimodularity.  That operator is the Dixmier-–Puk\' anszky 
Operator on $MAN$ and its domain is the Schwartz space $\cC(MAN)$ of
rapidly decreasing $C^\infty$ functions.

The kernel of the modular function $\delta$ of $MAN$ contains $MN$ and is
given on $A$ as follows.

\begin{lemma}\label{trace}
{\rm \cite[Lemmas 4.2 \& 4.3]{W2014}}
Let $\xi \in \ga$.  Then $\frac{1}{2}(\dim \gl_r + \dim \gz_r) \in \Z$
for $1\leqq r\leqq m$ and\\
\phantom{iiXX}{\rm (i)} the trace of $\ad(\xi)$ on $\gl_r$ is
  $\frac{1}{2}(\dim \gl_r + \dim \gz_r)\beta_r(\xi)$, \\
\phantom{iXX}{\rm (ii)} the trace of $\ad(\xi)$ on $\gn$ and on $\gp$ is
  $\frac{1}{2}\sum_r (\dim \gl_r + \dim \gz_r)\beta_r(\xi)$, and \\
\phantom{XX}{\rm (iii)}
  the determinant of $\Ad(\exp(\xi))$ on $\gn$ and on $\gp$ is
  $\prod_r \exp(\beta_r(\xi))^{\frac{1}{2} (\dim \gl_r + \dim \gz_r)}$.\\
The modular function $\delta = \Det\cdot\Ad: man \mapsto 
\prod_r \exp(\beta_r(\log a))^{\frac{1}{2} (\dim \gl_r + \dim \gz_r)}$.
\end{lemma}

Recall the {\em quasi-center determinant} 
$\Det_{\gs^*}(\lambda):=\prod_r (\beta_r(\lambda))^{\dim \gg_{\beta_r}}$\,.
It is a polynomial function on $\gs^*$, and (\cite[Proposition 4.7]{W2014})
the product $\Pf\cdot\Det_{\gs^*}$ is an $\Ad(MAN)$-semi-invariant
polynomial on $\gs^*$ of degree
$\frac{1}{2}(\dim \gn + \dim \gs)$ and of
weight equal to that of the modular function $\delta$.

Our fixed decomposition $\gn = \gv + \gs$ gives $N = VS$ where
$V = \exp(\gv)$ and $S = \exp(\gs)$.  Now define
\begin{equation}\label{defdp}
D: \text{ Fourier transform of } \Pf\cdot\Det_{\gs^*}
        \text{, acting on $MAN = MAVS$  by acting on $S$.}
\end{equation}
See \cite[Section 1]{W2014} for a discussion of the Schwartz space $\cC(MAN)$.

\begin{theorem}\label{dp-min-parab}{\rm (\cite[Theorem 4.9]{W2014})}
$D$ is an invertible self-adjoint differential operator of degree
$\frac{1}{2}(\dim \gn + \dim \gs)$ on $\cL^2(MAN)$ with dense
domain $\cC(MAN)$, and it is $\Ad(MAN)$-semi-invariant of
weight equal to the modular function $\delta$\,.
In other words $|D|$ is a Dixmier--Puk\' anszky Operator
on $MAN$ with domain $\cC(MAN)$.
\end{theorem}

For the Fourier Inversion Formula we also need to know the $\Ad^*(MA)$--orbits
on $\gt^*$.
\begin{proposition}\label{ma-orbits}
The $\Ad^*(MA)$--orbits on $\gt^*$ are the following.

In cases {\rm (1)} and {\rm (3)} of {\rm Theorem \ref{s-diagram-list}},
the number of $\Ad(MA)^*$--orbits on $\gt^*$ is $2^m$.   Fix nonzero
$\nu_r \in \gg_{\beta_r}$\,.  Then the orbits are the
$$
\cO_{(\varepsilon_1, \dots , \varepsilon_m)} =
\{\lambda = \lambda_1 + \dots + \lambda_m \mid 
\lambda_r \in \R^+\varepsilon_r\nu_r \text{ for } 1 \leqq r \leqq m
\text{ where each } \varepsilon_r = \pm 1\}.
$$
In case {\rm (2)} of {\rm Theorem \ref{s-diagram-list}}, there is just one
$\Ad(MA)^*$--orbit on $\gt^*$, i.e. $\Ad(MA)^*$ is transitive on $\gt^*$.
\end{proposition}
\noindent {\bf Proof.}  The assertions follow from Proposition \ref{m-action},
as follows.
In case (2) of Theorem \ref{s-diagram-list}, where $\dim \gg_{\beta_r} = 2$,
$\Ad^*(M)$ acts on $\gs^*$ by independent circle groups on the $\gg_{\beta_r}^*$
while $\Ad^*(A)$ acts on $\gs^*$ by independent positive scalar 
multiplication on the $\gg_{\beta_r}^*$\,.  In cases (1) and (3) of
Theorem \ref{s-diagram-list}, where $\dim \gg_{\beta_r} = 1$,
$\Ad^*(M)$ acts trivially on $\gs^*$ while $\Ad^*(A)$ acts on $\gs^*$ 
by independent positive scalar multiplication on the $\gg_{\beta_r}^*$\,.
\hfill $\square$

We combine Lemma \ref{ma-stabilizer}, Theorem \ref{sup-planch},
Theorem \ref{dp-min-parab} and Proposition \ref{ma-orbits} for the
Fourier Inversion Formula.  We need notation from Proposition \ref{ma-orbits}. 
For cases (1) and (3) of Theorem \ref{s-diagram-list} we fix nonzero 
$\nu_r \in \gg_{\beta_r}$; then for each
$\varepsilon = (\varepsilon_1, \dots , \varepsilon_m)$ we have
the orbit $\cO_\varepsilon$.  
As usual we write $\ell_x$ for left translate, $(\ell_xh)(y) = h(x^{-1}y)$ 
and $r_y$ for right translate $(r_yh)(x) = h(xy)$.
Using the structure of $M$ as a quotient of $\widetilde{M}$ by a
subgroup of the center of $\widetilde{G}$, from Proposition \ref{mtilde},
\cite[Theorem 6.1]{W2014} specializes as follows.

\begin{theorem}\label{fourier-inv}
Let $G$ be a real reductive Lie group whose minimal parabolic subgroup
$P = MAN$ is solvable.  Given a generic representation 
$\pi_{\chi,\alpha,\lambda} \in \widehat{MAN}$, its distribution character 
is tempered and is given by
$$
\Theta_{\pi_{\chi,\alpha,\lambda}}(f) = \tr \pi_{\chi,\alpha,\lambda}(f)
= \int_{M}\tr\chi(m) \int_{A} e^{i\alpha(\log a)} 
\Theta_{\pi_\lambda}(\ell_{(ma)^{-1}}f) da\, dm
\text{ for } f \in \cC(MAN).
$$
where $\Theta_{\pi_\lambda}$ is given by {\rm Theorem \ref{plancherel-general}}.
In Cases {\rm (1)} and {\rm (3)} of 
{\rm Theorem \ref{s-diagram-list}} the Fourier Inversion Formula is
$$
f(x) = c\int_{\chi \in \widehat{M}} 
	\left ( \int_{\alpha\in\ga_{\diamondsuit}^*} 
	\sum_\varepsilon 
	\left ( \int_{\lambda \in \cO_\varepsilon}
        \Theta_{\pi_{\chi,\alpha,\lambda}}(D(r(x)f)) |\Pf(\lambda)|d\lambda\right ) 
	d\alpha \right )\deg(\chi)\, d\chi
$$
where $c > 0$ depends on {\rm (\ref{c-d})(a)} and normalizations of Haar 
measures.  In Case {\rm (2)} of {\rm Theorem \ref{s-diagram-list}} the Fourier 
Inversion Formula is
$$
f(x) = c {\sum}_{\chi \in \widehat{M_\diamondsuit}}
	\left ( \int_{\alpha\in\ga_{\diamondsuit}^*}
	\left ( \int_{\lambda \in \gt^*}
        \Theta_{\pi_{\chi,\alpha,\lambda}}(D(r(x)f)) |\Pf(\lambda)|
	d\lambda \right ) d\alpha \right )
$$
where, again, $c > 0$ depends on {\rm (\ref{c-d})(a)} and normalizations of 
Haar measures.
\end{theorem}

\section{Appendix: Explicit Decompositions}\label{secA}
\setcounter{equation}{0}
One can see the decompositions of Section \ref{sec2} explicitly.  
This is closely related to the computations in \cite[Section 8]{W2012}. 
As noted in Corollary \ref{divisibility}, the only
case of Theorem \ref{s-diagram-list} where there is any divisibility is
$\gg = \gs\gu(\ell,\ell+1)$, and in that case the only divisibility is
given by the $\left \{ \tfrac{1}{2}\beta_r\,, \beta_r\right \} \subset 
\Delta^+(\gg,\ga)$.  For ease of terminology we say that $\Delta(\gg,\ga)$
is {\em nonmultipliable} if $\alpha \in \Delta(\gg,\ga)$ implies
$2\alpha \notin \Delta(\gg,\ga)$, {\em multipliable} otherwise.

We first consider the multipliable case 
$\gg = \gs\gu(\ell,\ell +1)$
\setlength{\unitlength}{.4 mm}
\begin{picture}(80,16)
\put(5,-1){\circle{2}}
\put(5,13){\circle{2}}
\put(4,4){$\updownarrow$}
\put(6,-1){\line(1,0){13}}
\put(6,13){\line(1,0){13}}
\put(22,-1){\circle*{1}}
\put(25,-1){\circle*{1}}
\put(28,-1){\circle*{1}}
\put(22,13){\circle*{1}}
\put(25,13){\circle*{1}}
\put(28,13){\circle*{1}}
\put(31,-1){\line(1,0){13}}
\put(31,13){\line(1,0){13}}
\put(45,-1){\circle{2}}
\put(45,13){\circle{2}}
\put(44,4){$\updownarrow$}
\put(46,-1){\line(1,0){6}}
\put(46,13){\line(1,0){6}}
\put(55,-1){\circle*{1}}
\put(58,-1){\circle*{1}}
\put(61,-1){\circle*{1}}
\put(55,13){\circle*{1}}
\put(58,13){\circle*{1}}
\put(61,13){\circle*{1}}
\put(64,-1){\line(1,0){6}}
\put(64,13){\line(1,0){6}}
\put(71,-1){\circle{2}}
\put(71,13){\circle{2}}
\put(70,4){$\updownarrow$}
\put(72,-1){\line(1,0){3}}
\put(72,13){\line(1,0){3}}
\put(75,-1){\line(0,1){14}}
\end{picture}
Here $\psi_i = \alpha_i|_{\ga}$ where $\{\alpha_1, \dots ,\alpha_{2\ell}\}$ are
the simple roots of $\Delta^+(\gg_{_\C},\gh_{_\C})$ in the usual order.
The multipliable roots in $\Delta^+(\gg,\ga)$ are just the
$\frac{1}{2}\beta_r = (\alpha_r + \dots + \alpha_\ell)|_\ga$\,, 
$1 \leqq r \leqq \ell$, where $\beta_r = 2(\psi_r + \dots + \psi_\ell) = 
(\alpha_r + \dots + \alpha_{2\ell -r+1})|_\ga$\,.  If $\alpha + \alpha' = 
\beta_r$ then, either $\alpha = \alpha' = \frac{1}{2}\beta_r$\,, or
one of $\alpha,\, \alpha'$ has form
$\gamma_{r,u} = \psi_r + \psi_{r+1} + \dots + \psi_u$ while the other is
$\gamma'_{r,u} = \psi_r + \psi_{r+1} + \dots + \psi_u + 2(\psi_{u+1} + \dots +
      \psi_\ell)$.
Now $\gl_r = \gg_{\beta_r} + \gg_{\frac{1}{2}\beta_r} +
\sum_{r \leqq u \leqq n}(\gg_{\gamma_{r,u}} + \gg_{\gamma'_{r,u}})$.
The conditions of (\ref{setup}) follow by inspection.

For the rest of this section we assume that $\gg \ne \gs\gu(\ell,\ell+1)$,
in other words, following Corollary \ref{divisibility}, that 
$\Delta(\gg,\ga)$ is nonmultipliable.

First suppose that $\Delta(\gg,\ga)$ is of type $A_{\ell - 1}$\,: 
\setlength{\unitlength}{.35 mm}
\begin{picture}(120,18)
\put(5,2){\circle{3}}
\put(2,5){$\psi_1$}
\put(6,2){\line(1,0){23}}
\put(30,2){\circle{3}}
\put(27,5){$\psi_2$}
\put(31,2){\line(1,0){23}}
\put(55,2){\circle{3}}
\put(52,5){$\psi_{\ell-2}$}
\put(56,2){\line(1,0){13}}
\put(74,2){\circle*{1}}
\put(77,2){\circle*{1}}
\put(80,2){\circle*{1}}
\put(86,2){\line(1,0){13}}
\put(100,2){\circle{3}}
\put(97,5){$\psi_{\ell - 1}$}
\end{picture}
Then $m = [\ell /2]$ and $\beta_r = \psi_r + \psi_{r+1} + \dots +
\psi_{\ell -r}$.  If $\alpha\,, \alpha' \in \Delta^+(\gg,\ga)$ with 
$\alpha + \alpha' = \beta_r$ then one of $\alpha\,, \alpha'$ must have form
$\gamma_{r,s} := \psi_r + \dots + \psi_s$ and the other must be 
$\gamma'_{r,s} := \psi_{s+1} + \dots + \psi_{\ell - r}$.  Thus
$\gl_r = \gg_{\beta_r} + \sum_{r \leqq s < \ell - r}\,
(\gg_{\gamma_{r,s}} + \gg_{\gamma'_{r,s}})$\,.
The conditions of (\ref{setup}) follow by inspection.

Next suppose that $\Delta(\gg,\ga)$ is of type $B_n$\,: 
\setlength{\unitlength}{.35 mm}
\begin{picture}(120,12)
\put(5,2){\circle{3}}
\put(2,5){$\psi_1$}
\put(6,2){\line(1,0){23}}
\put(30,2){\circle{3}}
\put(27,5){$\psi_2$}
\put(31,2){\line(1,0){13}}
\put(50,2){\circle*{1}}
\put(53,2){\circle*{1}}
\put(56,2){\circle*{1}}
\put(61,2){\line(1,0){13}}
\put(75,2){\circle{3}}
\put(72,5){$\psi_{n-1}$}
\put(76,2.5){\line(1,0){23}}
\put(76,1.5){\line(1,0){23}}
\put(100,2){\circle*{3}}
\put(97,5){$\psi_n$}
\end{picture}
Then $\beta_1 = \psi_1 + 2(\psi_2 + \dots + \psi_n)$, $\beta_2 = \psi_1$,
$\beta_3 = \psi_3 + 2(\psi_4 + \dots + \psi_n)$, $\beta_4 = \psi_3$, etc.
If $r$ is even, $\beta_r = \psi_{r-1}$ and $\gl_r = \gg_{\beta_r}$\,.

Now let $r$ be odd, $\beta_r = 
\psi_r + 2(\psi_{r+1} + \dots + \psi_n)$.  
If $\alpha,\, \alpha' \in \Delta^+(\gg,\ga)$ with $\alpha + \alpha'
= \beta_r$ then one possibility is that one of $\alpha\,, \alpha'$ has form
$\gamma_{r,u}:= \psi_r + \psi_{r+1} + \dots + \psi_u$ and the other is
$\gamma'_{r,u}:= \psi_{r+1} + \dots + \psi_u + 
		2(\psi_{u+1} + \dots + \psi_n)$ with $r < u < n$.
Another possibility is that one of $\alpha\,, \alpha'$ is 
$\gamma_{r,u} - \psi_r$ while the other is $\gamma'_{r,u} + \psi_r$\,.  
A third is that one of $\alpha\,, \alpha'$ is
$\gamma_{r,n}:= \psi_r + \psi_{r+1} + \dots + \psi_n$ while the other is
$\gamma'_{r,n}:= \psi_{r+1} + \dots + \psi_n$\,.  Then $\gl_r$ is the
sum of $\gg_{\beta_r}$ with the sum of all these possible 
$\gg_\alpha + \gg_{\alpha'}$, and the conditions of (\ref{setup}) follow
by inspection.

Let $\Delta(\gg,\ga)$ be of type $C_n$\,: 
\setlength{\unitlength}{.5 mm}
\begin{picture}(110,12)
\put(5,2){\circle*{3}}
\put(2,5){$\psi_1$}
\put(6,2){\line(1,0){23}}
\put(30,2){\circle*{3}}
\put(27,5){$\psi_2$}
\put(31,2){\line(1,0){13}}
\put(50,2){\circle*{1}}
\put(53,2){\circle*{1}}
\put(56,2){\circle*{1}}
\put(61,2){\line(1,0){13}}
\put(75,2){\circle*{3}}
\put(72,5){$\psi_{n-1}$}
\put(76,2.5){\line(1,0){23}}
\put(76,1.5){\line(1,0){23}}
\put(100,2){\circle{3}}
\put(97,5){$\psi_n$}
\end{picture}
Then $\beta_r = 2(\psi_r + \dots \psi_{n-1}) + \psi_n$ for $1 \leqq r < n$,
and $\beta_n = \psi_n$\,.  If $\alpha + \alpha' = \beta_r$ with $r < n$
then one of $\alpha,\, \alpha'$ has form 
$\gamma_{r,u} = \psi_r + \psi_{r+1} + \dots + \psi_u$ while the other is
$\gamma'_{r,u} = \psi_r + \psi_{r+1} + \dots + \psi_u + 2(\psi_{u+1} + \dots +
      \psi_{n-1}) + \psi_n$, $r \leqq u < n$. 
Note that $\gamma_{r,n-1} = \psi_r + \psi_{r+1} + \dots + \psi_{n-1}$ and
$\gamma'_{r,n-1} = \psi_r + \psi_{r+1} + \dots + \psi_n$\,.
Now $\gl_r = \gg_{\beta_r} + 
\sum_{r \leqq u < n}(\gg_{\gamma_{r,u}} + \gg_{\gamma'_{r,u}})$.
The conditions of (\ref{setup}) follow by inspection.

Let $\Delta(\gg,\ga)$ be of type $D_n\,$: 
\setlength{\unitlength}{.35 mm}
\begin{picture}(110,20)
\put(5,9){\circle{3}}
\put(2,12){$\psi_1$}
\put(6,9){\line(1,0){23}}
\put(30,9){\circle{3}}
\put(27,12){$\psi_2$}
\put(31,9){\line(1,0){13}}
\put(50,9){\circle*{1}}
\put(53,9){\circle*{1}}
\put(56,9){\circle*{1}}
\put(61,9){\line(1,0){13}}
\put(75,9){\circle{3}}
\put(70,12){$\psi_{n-2}$}
\put(76,8.5){\line(2,-1){13}}
\put(90,2){\circle{3}}
\put(93,0){$\psi_n$}
\put(76,9.5){\line(2,1){13}}
\put(90,16){\circle{3}}
\put(93,14){$\psi_{n-1}$}
\end{picture}
Then $\beta_1 = \psi_1 + 2(\psi_2 + \dots + \psi_{n-2}) + \psi_{n-1} + \psi_n$,
$\beta_2 = \psi_1$\,, 
$\beta_3 = \psi_3 + 2(\psi_4 + \dots + \psi_{n-2}) + \psi_{n-1} + \psi_n$,
$\beta_4 = \psi_3$\,, etc., until $r = n-3$.  

If $n$ is even then $m=n$,
$\beta_{n-3} = \psi_{n-3} + 2\psi_{n-2} + \psi_{n-1} + \psi_n$\,, 
$\beta_{n-2} = \psi_{n-3}$\,, 
$\beta_{n-1} = \psi_{n-1}$ and $\beta_n = \psi_n$\,.   
Then, if $r\leqq n-2$ is even we have $\beta_r = \psi_{r-1}$\,.  Thus
$\gl_r = \gg_{\beta_r}$ for $n$ even and either $r$ even or $r = n-1$. 

If $n$ is odd then
$m = n-1$, $\beta_{n-2} = \psi_{n-2} + \psi_{n-1} + \psi_n$ and 
$\beta_{n-1} = \psi_{n-2}$\,.  Thus $\beta_r = \psi_{r-1}$
and $\gl_r = \gg_{\beta_r}$ for $n$ odd and $r$ even.

That leaves the cases where $r$ is odd and $r \ne n-1$, so
$\beta_r = \psi_r + 2(\psi_{r+1} + \dots + \psi_{n-2}) + \psi_{n-1} + \psi_n$.
If $\alpha + \alpha' = \beta_r$\,, one possibility is that one of 
$\alpha,\, \alpha'$ is of the form
$\gamma_{r,u} := \psi_r + (\psi_{r+1} + \dots + \psi_u) \text{ with }
        r+1 \leqq u \leqq n-2$
while the other is $\gamma'_{r,u} := (\psi_{r+1} + \dots + \psi_u) +
  2(\psi_{u+1} + \dots + \psi_{n-2}) + \psi_{n-1} + \psi_n$\,, or that
one of $\alpha,\, \alpha'$ is of the form $\gamma_{r,u} - \psi_r$ while
the other is $\gamma'_{r,u} + \psi_r$\,.
A third possibility is that one of $\alpha,\, \alpha'$ is
$\gamma_{r,n-1} := \psi_r + \psi_{r+1} + \dots + \psi_{n-1}$ while
the other is $\gamma'_{r,n-1}:= \psi_{r+1} + \dots + \psi_{n-2} + \psi_n$\,.
The fourth possibility is that one of $\alpha,\, \alpha'$ is
$\gamma_{r,n-1} - \psi_{n-1} + \psi_n$ while the other is
$\gamma'_{r,n-1} + \psi_{n-1} - \psi_n$\,.  Then $\gl_r$ is the
sum of $\gg_{\beta_r}$ with the sum of all these possible
$\gg_\alpha + \gg_{\alpha'}$, and the conditions of (\ref{setup}) follow
by inspection.

Let $\Delta(\gg,\ga)$ be of type $G_2$ 
\setlength{\unitlength}{.75 mm}
\begin{picture}(25,10)
\put(5,3){\circle*{3}}
\put(3,-1){$\psi_1$}
\put(6,2){\line(1,0){13}}
\put(6,3){\line(1,0){13}}
\put(6,4){\line(1,0){13}}
\put(20,3){\circle{3}}
\put(18,-2){$\psi_2$}
\end{picture}.
Then $\beta_1 = 3\psi_1 + 2\psi_2$ and $\beta_2 = \psi_1$. If
$\alpha + \alpha' = \beta_1$ then either one of $\alpha,\, \alpha'$ is
$3\psi_1 + \psi_2$ and the other is $\psi_2$\,, or one of $\alpha,\, \alpha'$ is
$2\psi_1 + \psi_2$ and the other is $\psi_1 + \psi_2$\,.  Thus 
$\gl_1 = \gg_{\beta_1} + (\gg_{3\psi_1 + \psi_2} + \gg_{\psi_2}) +
(\gg_{2\psi_1 + \psi_2} + \gg_{\psi_1 + \Psi_2})$ and $\gl_2 = \gg_{\beta_2}$\,.
The conditions of (\ref{setup}) follow.

Suppose that $\Delta(\gg,\ga)$ is of type $F_4$
\setlength{\unitlength}{.75 mm}
\begin{picture}(60,13)
\put(10,3){\circle{3}}
\put(8,-2){$\psi_1$}
\put(11,3){\line(1,0){13}}
\put(25,3){\circle{3}}
\put(23,-2){$\psi_2$}
\put(26,2.5){\line(1,0){13}}
\put(26,3.5){\line(1,0){13}}
\put(40,3){\circle*{3}}
\put(38,-2){$\psi_3$}
\put(41,3){\line(1,0){13}}
\put(55,3){\circle*{3}}
\put(53,-2){$\psi_4$}
\end{picture}.
Then $\beta_1 = 2\psi_1 + 3\psi_2 + 4\psi_4 + 2\psi_4$,\, 
$\beta_2 = \psi_2 + 2\psi_3 + 2\psi_4$,\, $\beta_3 = \psi_2 + 2\psi_3$ and
$\beta_4 = \psi_2$. Thus 
$\gl_r = \gg_{\beta_r} + \sum_{(\gamma,\gamma') \in S_r}\,
(\gg_\gamma + \gg_{\gamma'})$ where
{\footnotesize
$$
\begin{aligned}
S_1 = \{&(\psi_1,\, \psi_1 + 3\psi_2 + 4\psi_3 + 2\psi_4),
	(\psi_1 + \psi_2,\, \psi_1 + 2\psi_2 + 4\psi_3 + 2\psi_4),\\
&(\psi_1 + \psi_2 + \psi_3,\, \psi_1 + 2\psi_2 + 3\psi_3 + 2\psi_4), 
(\psi_1 + \psi_2 + 2\psi_3,\,  \psi_1 + 2\psi_2 + 2\psi_3 + 2\psi_4), \\
&(\psi_1 + \psi_2 + \psi_3 + \psi_4,\, \psi_1 + 2\psi_2 + 3\psi_3 + \psi_4),
(\psi_1 + 2\psi_2 + 2\psi_3,\, \psi_1 + \psi_2 + 2\psi_3 + 2\psi_4),\\
&(\psi_1 + \psi_2 + 2\psi_3 + \psi_4,\,\psi_1 + 2\psi_2 + 2\psi_3 + \psi_4);\\
S_2 = \{&\{\psi_4,\, \psi_2 + 2\psi_3 + \psi_4\},  
      \{\psi_3 + \psi_4,\, \psi_2 + \psi_3 + \psi_4\}\};\,
S_3 = \{\psi_3,\, \psi_2 + \psi_3\}\, \text{ and } S_4 = \emptyset
\end{aligned}
$$
}
The conditions of (\ref{setup}) follow.

Suppose that $\Delta(\gg,\ga)$ is of type $E_6$ 
\setlength{\unitlength}{.35 mm}
\begin{picture}(70,23)
\put(10,12){\circle{3}}
\put(8,15){$\psi_1$}
\put(11,12){\line(1,0){13}}
\put(25,12){\circle{3}}
\put(23,15){$\psi_3$}
\put(26,12){\line(1,0){13}}
\put(40,12){\circle{3}}
\put(38,15){$\psi_4$}
\put(41,12){\line(1,0){13}}
\put(55,12){\circle{3}}
\put(53,15){$\psi_5$}
\put(56,12){\line(1,0){13}}
\put(70,12){\circle{3}}
\put(68,15){$\psi_6$}
\put(40,11){\line(0,-1){13}}
\put(40,-3){\circle{3}}
\put(43,-3){$\psi_2$}
\end{picture}
Then the strongly orthogonal roots $\beta_i$ are given by 
$\beta_1 = \psi_1 + 2\psi_2 + 2\psi_3 + 3\psi_4 + 2\psi_5 + \psi_6$\,,
$\beta_2 = \psi_1 + \psi_3 + \psi_4 + \psi_5 + \psi_6,$\,,
$\beta_3 = \psi_3 + \psi_4 + \psi_5$  and  
$\beta_4 = \psi_4$\,.  
Now $\gl_r = \gg_{\beta_r} + \sum_{(\gamma,\gamma') \in S_r} 
(\gg_\gamma + \gg_{\gamma'})$ where
{\footnotesize
$$
\begin{aligned}
S_1 = \{&(\psi_2, \psi_1 + \psi_2 + 2\psi_2 + 3\psi_4 + 2\psi_5 + \psi_6),
         (\psi_2 + \psi_4, \psi_1 + \psi_2 + 2\psi_3 + 2\psi_4 + 2\psi_5 
		+ \psi_6),\\
         &(\psi_2 + \psi_3 + \psi_4, \psi_1 + \psi_2 + \psi_3 + 2\psi_4 + 
		2\psi_5 + \psi_6),
         (\psi_2 + \psi_4 + \psi_5, \psi_1 + \psi_2 + 2\psi_3 + 2\psi_4 + 
		\psi_5 + \psi_6),\\
         &(\psi_1 + \psi_2 + \psi_3 + \psi_4, \psi_2 + \psi_3 + 2\psi_4 + 
		2\psi_5 + \psi_6),
         (\psi_2 + \psi_3 + \psi_4 + \psi_5, \psi_1 + \psi_2 + \psi_3 + 
		2\psi_4 + \psi_5 + \psi_6),\\
         &(\psi_2 + \psi_4 + \psi_5 + \psi_6, \psi_1 + \psi_2 + 2\psi_3 + 
		2\psi_4 + \psi_5),
         (\psi_1 + \psi_2 + \psi_3 + \psi_4 + \psi_5, \psi_2 + \psi_3 + 
		2\psi_4 + \psi_5 + \psi_6),\\
         &(\psi_2 + \psi_3 + 2\psi_4 + \psi_5, \psi_1 + \psi_2 + \psi_3 + 
		\psi_4 + \psi_5 + \psi_6),
         (\psi_2 + \psi_3 + \psi_4 + \psi_5 + \psi_6,  \psi_1 + \psi_2 + 
		\psi_3 + 2\psi_4 + \psi_5) \};\\
S_2 = \{ &(\psi_1, \psi_3 + \psi_4 + \psi_5 + \psi_6,
          \psi_6, \psi_1 + \psi_3 + \psi_4 + \psi_5), 
          (\psi_1 + \psi_3, \psi_4 + \psi_5 + \psi_6,
          \psi_5 + \psi_6, \psi_1 + \psi_3 + \psi_4) \};\\
S_3 = \{&(\psi_3, \psi_4 + \psi_5), (\psi_5, \psi_3 + \psi_4)\}
\text{ and } S_4 = \emptyset.
\end{aligned}
$$
}
The conditions of (\ref{setup}) follow.

Next, suppose that $\Delta(\gg,\ga)$ is of type $E_7$  
\setlength{\unitlength}{.35 mm}
\begin{picture}(100,23)
\put(10,12){\circle{3}}
\put(8,15){$\psi_1$}
\put(11,12){\line(1,0){13}}
\put(25,12){\circle{3}}
\put(23,15){$\psi_3$}
\put(26,12){\line(1,0){13}}
\put(40,12){\circle{3}}
\put(38,15){$\psi_4$}
\put(41,12){\line(1,0){13}}
\put(55,12){\circle{3}}
\put(53,15){$\psi_5$}
\put(56,12){\line(1,0){13}}
\put(70,12){\circle{3}}
\put(68,15){$\psi_6$}
\put(71,12){\line(1,0){13}}
\put(85,12){\circle{3}}
\put(83,15){$\psi_7$}
\put(40,11){\line(0,-1){13}}
\put(40,-3){\circle{3}}
\put(43,-3){$\psi_2$}
\end{picture}.
Then the strongly orthogonal roots are
$\beta_1 = 2\psi_1 + 2\psi_2 + 3\psi_3+4\psi_4+3\psi_5 + 2\psi_6 + \psi_7$\,,
$\beta_2 = \psi_2 + \psi_3 + 2\psi_4 + 2\psi_5 + 2\psi_6 + \psi_7$\,,
$\beta_3 = \psi_7$\,,
$\beta_4 = \psi_2 + \psi_3 + 2\psi_4 + \psi_5$\,,
$\beta_5 = \psi_2$\,,
$\beta_6 = \psi_3$ and 
$\beta_7 = \psi_5$\,.  
Now $\gl_r = \gg_{\beta_r} + \sum_{(\gamma,\gamma') \in S_r}
(\gg_\gamma + \gg_{\gamma'})$ where
\footnotesize
$$
\begin{aligned}
S_1 = \{&( \psi_1,\, 
	\psi_1 + 2\psi_2 + 3\psi_3 + 4\psi_4 + 3\psi_5 + 2\psi_6 + \psi_7),\\
        &( \psi_1 + \psi_3,\,
	\psi_1 + 2\psi_2 + 2\psi_3 + 4\psi_4 + 3\psi_5 + 2\psi_6 + \psi_7),\\
        &( \psi_1 + \psi_3 + \psi_4,\,
	 \psi_1 + 2\psi_2 + 2\psi_3 + 3\psi_4 + 3\psi_5 + 2\psi_6 + \psi_7),\\
        &( \psi_1 + \psi_2 + \psi_3 + \psi_4,\,
	 \psi_1 + \psi_2 + 2\psi_3 + 3\psi_4 + 3\psi_5 + 2\psi_6 + \psi_7),\\
        &( \psi_1 + \psi_3 + \psi_4 + \psi_5,\,
	 \psi_1 + 2\psi_2 + 2\psi_3 + 3\psi_4 + 2\psi_5 + 2\psi_6 + \psi_7),\\
        &( \psi_1 + \psi_2 + \psi_3 + \psi_4 + \psi_5,\,
	 \psi_1 + \psi_2 + 2\psi_3 + 3\psi_4 + 2\psi_5 + 2\psi_6 + \psi_7),\\
        &( \psi_1 + \psi_3 + \psi_4 + \psi_5 + \psi_6,\,
	 \psi_1 + 2\psi_2 + 2\psi_3 + 3\psi_4 + 2\psi_5 + \psi_6 + \psi_7),\\
        &( \psi_1 + \psi_2 + \psi_3 + 2\psi_4 + \psi_5,\,
	 \psi_1 + \psi_2 + 2\psi_3 + 2\psi_4 + 2\psi_5 + 2\psi_6 + \psi_7),\\
        &( \psi_1 + \psi_2 + \psi_3 + \psi_4 + \psi_5 + \psi_6,\,
	 \psi_1 + \psi_2 + 2\psi_3 + 3\psi_4 + 2\psi_5 + \psi_6 + \psi_7),\\
        &( \psi_1 + \psi_3 + \psi_4 + \psi_5 + \psi_6 + \psi_7,\,
	 \psi_1 + 2\psi_2 + 2\psi_3 + 3\psi_4 + 2\psi_5 + \psi_6),\\
        &( \psi_1 + \psi_2 + 2\psi_3 + 2\psi_4 + \psi_5,\,
	 \psi_1 + \psi_2 + \psi_3 + 2\psi_4 + 2\psi_5 + 2\psi_6 + \psi_7),\\
        &( \psi_1 + \psi_2 + \psi_3 + 2\psi_4 + \psi_5 + \psi_6,\,
	 \psi_1 + \psi_2 + 2\psi_3 + 2\psi_4 + 2\psi_5 + \psi_6 + \psi_7),\\
        &( \psi_1 + \psi_2 + \psi_3 + \psi_4 + \psi_5 + \psi_6 + \psi_7,\,
	 \psi_1 + \psi_2 + 2\psi_3 + 3\psi_4 + 2\psi_5 + \psi_6),\\
        &( \psi_1 + \psi_2 + 2\psi_3 + 2\psi_4 + \psi_5 + \psi_6,\,
	 \psi_1 + \psi_2 + \psi_3 + 2\psi_4 + 2\psi_5 + \psi_6 + \psi_7),\\
        &( \psi_1 + \psi_2 + \psi_3 + 2\psi_4 + 2\psi_5 + \psi_6,\,
	 \psi_1 + \psi_2 + 2\psi_3 + 2\psi_4 + \psi_5 + \psi_6 + \psi_7),\\
        &( \psi_1 + \psi_2 + \psi_3 + 2\psi_4 + \psi_5 + \psi_6 + \psi_7,\,
	 \psi_1 + \psi_2 + 2\psi_3 + 2\psi_4 + 2\psi_5 + \psi_6)\}
\end{aligned}
$$
while
{\footnotesize
$$
\begin{aligned}
S_2 = \{&(\psi_6,
	     \psi_2 + \psi_3 + 2\psi_4 + 2\psi_5 + \psi_6 + \psi_7),
        ( \psi_5 + \psi_6,
	     \psi_2 + \psi_3 + 2\psi_4 + \psi_5 + \psi_6 + \psi_7),\\
        &(\psi_6 + \psi_7,
	     \psi_2 + \psi_3 + 2\psi_4 + 2\psi_5 + \psi_6),
        ( \psi_4 + \psi_5 + \psi_6,
	     \psi_2 + \psi_3 + \psi_4 + \psi_5 + \psi_6 + \psi_7), \\
        &(\psi_5 + \psi_6 + \psi_7,
	     \psi_2 + \psi_3 + 2\psi_4 + \psi_5 + \psi_6),
        ( \psi_2 + \psi_4 + \psi_5 + \psi_6,
	     \psi_3 + \psi_4 + \psi_5 + \psi_6 + \psi_7),\\
        &(\psi_3 + \psi_4 + \psi_5 + \psi_6,
	     \psi_2 + \psi_4 + \psi_5 + \psi_6 + \psi_7),
        ( \psi_4 + \psi_5 + \psi_6 + \psi_7,
	     \psi_2 + \psi_3 + \psi_4 + \psi_5 + \psi_6) \}\\
S_4 = \{&(\psi_4, \,\psi_2 + \psi_3 + \psi_4 + \psi_5),
        (\psi_2 + \psi_4, \,\psi_3 + \psi_4 + \psi_5),
        (\psi_3 + \psi_4, \,\psi_2 + \psi_4 + \psi_5),
        (\psi_4 + \psi_5, \,\psi_2 + \psi_3 + \psi_4)\}
\end{aligned}
$$
}
and $S_3 = S_5 = S_6 = S_7 = \emptyset$.
The conditions of (\ref{setup}) follow.

Finally suppose that $\Delta(\gg,\ga)$ is of type $E_8$
\setlength{\unitlength}{.35 mm}
\begin{picture}(110,23)
\put(10,12){\circle{3}}
\put(8,15){$\psi_1$}
\put(11,12){\line(1,0){13}}
\put(25,12){\circle{3}}
\put(23,15){$\psi_3$}
\put(26,12){\line(1,0){13}}
\put(40,12){\circle{3}}
\put(38,15){$\psi_4$}
\put(41,12){\line(1,0){13}}
\put(55,12){\circle{3}}
\put(53,15){$\psi_5$}
\put(56,12){\line(1,0){13}}
\put(70,12){\circle{3}}
\put(68,15){$\psi_6$}
\put(71,12){\line(1,0){13}}
\put(85,12){\circle{3}}
\put(83,15){$\psi_7$}
\put(86,12){\line(1,0){13}}
\put(100,12){\circle{3}}
\put(98,15){$\psi_8$}
\put(40,11){\line(0,-1){13}}
\put(40,-3){\circle{3}}
\put(43,-3){$\psi_2$}
\end{picture}.
Then 
$\beta_1 = 2\psi_1 + 3\psi_2 + 4\psi_3 + 6\psi_4 + 5\psi_5
	+ 4\psi_6 + 3\psi_7 + 2\psi_8$\,,
$\beta_2 = 2\psi_1 + 2\psi_2 + 3\psi_3 + 4\psi_4 + 3\psi_5 +2\psi_6 + \psi_7$\,,
$\beta_3 = \psi_2 + \psi_3 + 2\psi_4 + 2\psi_5 + 2\psi_6 + \psi_7$\,,
$\beta_4 = \psi_7$\,,
$\beta_5 = \psi_2 + \psi_3 + 2\psi_4 + \psi_5$\,,
$\beta_6 = \psi_2$\,,
$\beta_7 = \psi_3$\,,  and 
$\beta_8 = \psi_5$\,.
Now $\gl_r = \gg_{\beta_r} + \sum_{(\gamma,\gamma') \in S_r} 
	(\gg_\gamma + \gg_{\gamma'})$ where
{\footnotesize
$$
\begin{aligned}
S_4 = S_6 &= S_7 = S_8 = \emptyset \\
S_5 = \{&(\psi_4,\,           \psi_2+\psi_3+\psi_4+\psi_5),
        (\psi_2+\psi_4,\,    \psi_3+\psi_4+\psi_5), \\
        &(\psi_3+\psi_4,\,    \psi_2+\psi_4+\psi_5),
        (\psi_4+\psi_5,\,    \psi_2+\psi_3+\psi_4 )\};\\
S_3 = \{&(\psi_6,\,   \psi_2+\psi_3+2\psi_4+2\psi_5+\psi_6+\psi_7 ),
  (\psi_5+\psi_6 ,\,  \psi_2+\psi_3+2\psi_4+\psi_5+\psi_6+\psi_7 ), \\
  &(\psi_6+\psi_7 ,\,  \psi_2+\psi_3+2\psi_4+2\psi_5+\psi_6 ),
  (\psi_4+\psi_5+\psi_6 ,\,  \psi_2+\psi_3+\psi_4+\psi_5+\psi_6+\psi_7 ), \\
  &(\psi_5+\psi_6+\psi_7 ,\,  \psi_2+\psi_3+2\psi_4+\psi_5+\psi_6 ),
  (\psi_2+\psi_4+\psi_5+\psi_6 ,\,  \psi_3+\psi_4+\psi_5+\psi_6+\psi_7 ), \\
  &(\psi_3+\psi_4+\psi_5+\psi_6 ,\,  \psi_2+\psi_4+\psi_5+\psi_6+\psi_7 ),
  (\psi_4+\psi_5+\psi_6+\psi_7 ,\,  \psi_2+\psi_3+\psi_4+\psi_5+\psi_6 )\};
\end{aligned}
$$

$$
\begin{aligned}
S_1 = \{&(\psi_8 ,\phantom{X}
		2\psi_1+3\psi_2+4\psi_3+6\psi_4+5\psi_5+4\psi_6+3\psi_7+\psi_8),\\
        &(\psi_7+\psi_8 ,\phantom{X}
		2\psi_1+3\psi_2+4\psi_3+6\psi_4+5\psi_5+4\psi_6+2\psi_7+\psi_8),\\
        &(\psi_6+\psi_7+\psi_8 ,\phantom{X}
		2\psi_1+3\psi_2+4\psi_3+6\psi_4+5\psi_5+3\psi_6+2\psi_7+\psi_8),\\
	&(\psi_5+\psi_6+\psi_7+\psi_8 ,\phantom{X}
		2\psi_1+3\psi_2+4\psi_3+6\psi_4+4\psi_5+3\psi_6+2\psi_7+\psi_8),\\
	&(\psi_4+\psi_5+\psi_6+\psi_7+\psi_8 ,\phantom{X}
		2\psi_1+3\psi_2+4\psi_3+5\psi_4+4\psi_5+3\psi_6+2\psi_7+\psi_8),\\
        &(\psi_2+\psi_4+\psi_5+\psi_6+\psi_7+\psi_8 ,\phantom{X}
		2\psi_1+2\psi_2+4\psi_3+5\psi_4+4\psi_5+3\psi_6+2\psi_7+\psi_8),\\
        &(\psi_3+\psi_4+\psi_5+\psi_6+\psi_7+\psi_8 ,\phantom{X}
		2\psi_1+3\psi_2+3\psi_3+5\psi_4+4\psi_5+3\psi_6+2\psi_7+\psi_8),\\
        &(\psi_1+\psi_3+\psi_4+\psi_5+\psi_6+\psi_7+\psi_8 ,\phantom{X}
		\psi_1+3\psi_2+3\psi_3+5\psi_4+4\psi_5+3\psi_6+2\psi_7+\psi_8),\\
        &(\psi_2+\psi_3+\psi_4+\psi_5+\psi_6+\psi_7+\psi_8 ,\phantom{X}
		2\psi_1+2\psi_2+3\psi_3+5\psi_4+4\psi_5+3\psi_6+2\psi_7+\psi_8),\\
        &(\psi_1+\psi_2+\psi_3+\psi_4+\psi_5+\psi_6+\psi_7+\psi_8 ,\phantom{X}
		\psi_1+2\psi_2+3\psi_3+5\psi_4+4\psi_5+3\psi_6+2\psi_7+\psi_8),\\
	&(\psi_2+\psi_3+2\psi_4+\psi_5+\psi_6+\psi_7+\psi_8 ,\phantom{X}
		2\psi_1+2\psi_2+3\psi_3+4\psi_4+4\psi_5+3\psi_6+2\psi_7+\psi_8),\\
	&(\psi_1+\psi_2+\psi_3+2\psi_4+\psi_5+\psi_6+\psi_7+\psi_8 ,\phantom{X}
		\psi_1+2\psi_2+3\psi_3+4\psi_4+4\psi_5+3\psi_6+2\psi_7+\psi_8),\\
	&(\psi_2+\psi_3+2\psi_4+2\psi_5+\psi_6+\psi_7+\psi_8 ,\phantom{X}
		2\psi_1+2\psi_2+3\psi_3+4\psi_4+3\psi_5+3\psi_6+2\psi_7+\psi_8),\\
	&(\psi_1+\psi_2+2\psi_3+2\psi_4+\psi_5+\psi_6+\psi_7+\psi_8 ,\phantom{X}
		\psi_1+2\psi_2+2\psi_3+4\psi_4+4\psi_5+3\psi_6+2\psi_7+\psi_8),\\
	&(\psi_1+\psi_2+\psi_3+2\psi_4+2\psi_5+\psi_6+\psi_7+\psi_8 ,\phantom{X}
		\psi_1+2\psi_2+3\psi_3+4\psi_4+3\psi_5+3\psi_6+2\psi_7+\psi_8),\\
	&(\psi_2+\psi_3+2\psi_4+2\psi_5+2\psi_6+\psi_7+\psi_8 ,\phantom{X}
		2\psi_1+2\psi_2+3\psi_3+4\psi_4+3\psi_5+2\psi_6+2\psi_7+\psi_8),\\
	&(\psi_1+\psi_2+2\psi_3+2\psi_4+2\psi_5+\psi_6+\psi_7+\psi_8 ,\phantom{X}
		\psi_1+2\psi_2+2\psi_3+4\psi_4+3\psi_5+3\psi_6+2\psi_7+\psi_8),\\
	&(\psi_1+\psi_2+\psi_3+2\psi_4+2\psi_5+2\psi_6+\psi_7+\psi_8 ,\phantom{X}
		\psi_1+2\psi_2+3\psi_3+4\psi_4+3\psi_5+2\psi_6+2\psi_7+\psi_8),\\
	&(\psi_2+\psi_3+2\psi_4+2\psi_5+2\psi_6+2\psi_7+\psi_8 ,\phantom{X}
		2\psi_1+2\psi_2+3\psi_3+4\psi_4+3\psi_5+2\psi_6+\psi_7+\psi_8),\\
	&( \psi_1+\psi_2+2\psi_3+3\psi_4+2\psi_5+\psi_6+\psi_7+\psi_8 ,\phantom{X}
		\psi_1+2\psi_2+2\psi_3+3\psi_4+3\psi_5+3\psi_6+2\psi_7+\psi_8),\\
	&( \psi_1+\psi_2+2\psi_3+2\psi_4+2\psi_5+2\psi_6+\psi_7+\psi_8 ,\phantom{X}
		\psi_1+2\psi_2+2\psi_3+4\psi_4+3\psi_5+2\psi_6+2\psi_7+\psi_8),\\
	&( \psi_1+\psi_2+\psi_3+2\psi_4+2\psi_5+2\psi_6+2\psi_7+\psi_8 ,\phantom{X}
		\psi_1+2\psi_2+3\psi_3+4\psi_4+3\psi_5+2\psi_6+\psi_7+\psi_8),\\
	&( \psi_1+2\psi_2+2\psi_3+3\psi_4+2\psi_5+\psi_6+\psi_7+\psi_8 ,\phantom{X}
		\psi_1+\psi_2+2\psi_3+3\psi_4+3\psi_5+3\psi_6+2\psi_7+\psi_8),\\
	&( \psi_1+\psi_2+2\psi_3+3\psi_4+2\psi_5+2\psi_6+\psi_7+\psi_8 ,\phantom{X}
		\psi_1+2\psi_2+2\psi_3+3\psi_4+3\psi_5+2\psi_6+2\psi_7+\psi_8),\\
	&( \psi_1+\psi_2+2\psi_3+2\psi_4+2\psi_5+2\psi_6+2\psi_7+\psi_8 ,\phantom{X}
		\psi_1+2\psi_2+2\psi_3+4\psi_4+3\psi_5+2\psi_6+\psi_7+\psi_8),\\
	&( \psi_1+2\psi_2+2\psi_3+3\psi_4+2\psi_5+2\psi_6+\psi_7+\psi_8 ,\phantom{X}
		\psi_1+\psi_2+2\psi_3+3\psi_4+3\psi_5+2\psi_6+2\psi_7+\psi_8),\\
	&( \psi_1+\psi_2+2\psi_3+3\psi_4+3\psi_5+2\psi_6+\psi_7+\psi_8 ,\phantom{X}
		\psi_1+2\psi_2+2\psi_3+3\psi_4+2\psi_5+2\psi_6+2\psi_7+\psi_8),\\
	&( \psi_1+\psi_2+2\psi_3+3\psi_4+2\psi_5+2\psi_6+2\psi_7+\psi_8 ,\phantom{X}
		\psi_1+2\psi_2+2\psi_3+3\psi_4+3\psi_5+2\psi_6+\psi_7+\psi_8 )\};
\end{aligned}
$$
}
while
{\footnotesize
$$
\begin{aligned}
S_2 = \{&(\psi_1, \phantom{X}
		\psi_1+2\psi_2+3\psi_3+4\psi_4+3\psi_5+2\psi_6+\psi_7),\\
	&(\psi_1+\psi_3, \phantom{X}
		\psi_1+2\psi_2+2\psi_3+4\psi_4+3\psi_5+2\psi_6+\psi_7),\\
	&(\psi_1+\psi_3+\psi_4, \phantom{X}
		\psi_1+2\psi_2+2\psi_3+3\psi_4+3\psi_5+2\psi_6+\psi_7),\\
	&(\psi_1+\psi_2+\psi_3+\psi_4, \phantom{X}
		\psi_1+\psi_2+2\psi_3+3\psi_4+3\psi_5+2\psi_6+\psi_7),\\
	&(\psi_1+\psi_3+\psi_4+\psi_5, \phantom{X}
		\psi_1+2\psi_2+2\psi_3+3\psi_4+2\psi_5+2\psi_6+\psi_7),\\
	&(\psi_1+\psi_2+\psi_3+\psi_4+\psi_5, \phantom{X}
		\psi_1+\psi_2+2\psi_3+3\psi_4+2\psi_5+2\psi_6+\psi_7),\\
	&(\psi_1+\psi_3+\psi_4+\psi_5+\psi_6, \phantom{X}
		\psi_1+2\psi_2+2\psi_3+3\psi_4+2\psi_5+\psi_6+\psi_7),\\
	&(\psi_1+\psi_2+\psi_3+2\psi_4+\psi_5, \phantom{X}
		\psi_1+\psi_2+2\psi_3+2\psi_4+2\psi_5+2\psi_6+\psi_7),\\
	&(\psi_1+\psi_2+\psi_3+\psi_4+\psi_5+\psi_6, \phantom{X}
		\psi_1+\psi_2+2\psi_3+3\psi_4+2\psi_5+\psi_6+\psi_7),\\
	&(\psi_1+\psi_3+\psi_4+\psi_5+\psi_6+\psi_7, \phantom{X}
		\psi_1+2\psi_2+2\psi_3+3\psi_4+2\psi_5+\psi_6),\\
	&(\psi_1+\psi_2+2\psi_3+2\psi_4+\psi_5, \phantom{X}
		\psi_1+\psi_2+\psi_3+2\psi_4+2\psi_5+2\psi_6+\psi_7),\\
	&(\psi_1+\psi_2+\psi_3+2\psi_4+\psi_5+\psi_6, \phantom{X}
		\psi_1+\psi_2+2\psi_3+2\psi_4+2\psi_5+\psi_6+\psi_7),\\
	&(\psi_1+\psi_2+\psi_3+\psi_4+\psi_5+\psi_6+\psi_7, \phantom{X}
		\psi_1+\psi_2+2\psi_3+3\psi_4+2\psi_5+\psi_6),\\
	&(\psi_1+\psi_2+2\psi_3+2\psi_4+\psi_5+\psi_6, \phantom{X}
		\psi_1+\psi_2+\psi_3+2\psi_4+2\psi_5+\psi_6+\psi_7),\\
	&(\psi_1+\psi_2+\psi_3+2\psi_4+2\psi_5+\psi_6, \phantom{X}
		\psi_1+\psi_2+2\psi_3+2\psi_4+\psi_5+\psi_6+\psi_7),\\
	&(\psi_1+\psi_2+\psi_3+2\psi_4+\psi_5+\psi_6+\psi_7, \phantom{X}
		\psi_1+\psi_2+2\psi_3+2\psi_4+2\psi_5+\psi_6 )\}.
\end{aligned}
$$
}

\enddocument
\end